\documentclass[preprint,number,sort&compress]{elsarticle}
\usepackage{amsmath,amsfonts}
\usepackage{hyperref}
\usepackage{rotating}
\usepackage{lscape}
\newcommand{\RR}{\mathbb R}
\newcommand{\CC}{\mathbb C}
\newcommand{\Fra}[2]{\displaystyle \frac{ #1}{ #2}}
\newcommand{\rmi}{\mathrm{i}}
\newcommand{\rme}{\mathrm{e}}
\newcommand{\de}{\,\mathrm{d}}
\newcommand{\abs}[1]{\left\lvert#1\right\rvert}

\newcommand{\Grad}[1]{\displaystyle \boldsymbol{\nabla} #1}
\newcommand{\Lapl}[1]{\displaystyle {\nabla}^2 #1}
\newcommand{\Div}[1]{\displaystyle \boldsymbol{\nabla} \cdot #1}

\newcommand{\Vmy}[1]{\mathbf{#1}}
\newcommand{\Dpar}[2]{\displaystyle \Fra{\partial #1}{\partial #2}}
\title{Reliability of the time splitting Fourier method\\
for singular solutions in quantum fluids}
\author[vr]{M.~Caliari\corref{cor1}}
\ead{marco.caliari@univr.it}
\author[vr]{S.~Zuccher}
\ead{simone.zuccher@univr.it}
\cortext[cor1]{Corresponding author}
\address[vr]{Department of Computer Science, University of Verona,
Strada Le Grazie 15, 37134~Verona, Italy}
\begin{document}
\begin{abstract}
We extensively study the numerical accuracy of the well-known time splitting 
Fourier spectral method for the approximation of singular solutions of the
Gross--Pitaevskii equation. In particular, we explore its capability of
preserving a steady-state vortex solution, whose density profile
is approximated by a very accurate diagonal Pad\'e expansion
of order 8, here explicitly derived for the first time.
Although the Fourier spectral method turns out to be only slightly more 
accurate than a time splitting finite difference scheme, the former is
reliable and efficient. 
Moreover, at a post-processing stage, it allows an accurate evaluation of 
the solution outside grid points, thus becoming particularly appealing when
high resolution is needed, such as in the study of quantum vortex interactions.
\end{abstract}
\maketitle
\section{Introduction}
Quantum turbulence~\cite{VIN2008,PL2011,BSS2014}, as well as classical
turbulence~\cite{F1995,P2000}, is dominated by reconnection of 
vortical structures which is much simpler to treat in the framework of quantum 
fluids rather than in viscous fluids~\cite{HD2011}, while leading to similar 
features such as time asymmetry~\cite{ZCBB12}.
Despite the fundamental differences between the two forms of turbulence,
there are reasons to believe that the understanding of quantum turbulence 
might shed new light on the understanding of its classical
counterpart~\cite{BSS2014}. 

Quantum fluids dynamics is properly described by the
Gross--Pitaevskii equation~(GPE)~\cite{P1961,G1963}
\begin{equation}\label{eq:GPE}
\frac{\partial \psi}{\partial t} = 
\Fra \rmi 2 \Lapl \psi + \Fra \rmi 2\left(1 - \abs{\psi}^2\right)\psi,
\end{equation}%
where $\psi$ is the complex wave function.
Through the Madelung transformation 
$\psi = \sqrt\rho\, \mathrm{exp}(\rmi \theta)$, equation~\eqref{eq:GPE} 
can be viewed in classical fluid dynamical terms as
\begin{eqnarray}
\label{eq:GPEinNScont}
  \Dpar{\rho}{t}+\Dpar{(\rho u_j)}{x_j}&=&0,\\
\label{eq:GPEinNSmom}
  \rho \left(\Dpar{u_i}{t}+u_j\Dpar{u_i}{u_j} \right)& = &
  -\Dpar{p}{x_i} + \Dpar{\tau_{ij}}{x_j},
\end{eqnarray}
where $\rho=|\psi|^2$ denotes density, ${\boldsymbol u}=\Grad{\theta}$ 
velocity, $p=\frac{\rho^2}{4}$ pressure, and
$\tau_{ij}=\frac 1 4 \rho\frac{\partial^2\ln\rho}{\partial x_i \partial x_j}$ 
the so-called  quantum stress ($i,j=1,2,3$). 
Defects in the wave function $\psi$ are interpreted as infinitesimally thin 
vortices of constant circulation 
$\Gamma=\oint{\boldsymbol u}\cdot\!\de{\boldsymbol s}=2\pi$, with  healing
length $\xi=1$.  
GPE conserves the mass and the energy
\begin{equation}\label{eq:energy}
E=\frac{1}{2}\int \abs{\nabla \psi}^2\de x+
\frac{1}{4}\int(1-\abs{\psi}^2)^2\de x.
\end{equation}

The main reason for preferring the GPE approach to others for the study of
quantum turbulence is that it guarantees
a natural dynamics of interacting vortices~\cite{ZR15} while resolving fine
scales up to the vortex core~\cite{BSS2014,KLP2014}.
On the contrary, methods based on the inviscid Euler equations (either 
their direct numerical simulation~\cite{BK08} or 
vortex filament methods~\cite{HB2014}) are unable to automatically perform
vortex reconnections, being forbidden by Euler dynamics.

The numerical solution of the GPE~\eqref{eq:GPE} 
is normally carried out by employing
time splitting Fourier methods~\cite{KL93,ZCBB12,AZCPPB14,ZR15} and by imposing
vortices in the form of singular phase defects in a unitary background density,
i.e. $\rho(x)=\abs{\psi(x)}^2\to 1$ when $\abs{x}\to \infty$.
However, these methods assume periodic boundary conditions.
Solutions which are not periodic must be mirrored in the directions 
lacking periodicity~\cite{KL93}, thus imposing doubling of the degrees of 
freedom in each of those directions and a consequent increase of the
computational effort.

Recent developments~\cite{SKP14,ZR15} have shown that reconnections in 
quantum fluids
are strictly related to topological features characterizing the interacting
vortex tubes such as writhe, total torsion and intrinsic twist.
These quantities depend on the fine details of the curve that describes the 
vortex centerline (its third derivative with respect to the curvilinear 
abscissa is required for computing torsion) and on the phase of the wavefunction
$\psi$ in the neighborhood of the vortex centerline.
Therefore 
it is 
paramount to resort to high resolution numerical simulations of 
equation~\eqref{eq:GPE}, especially in the proximity of the 
reconnection event.

With the goal of assessing the goodness of time-splitting Fourier methods
for singular solutions on uniform grids versus time-splitting finite
differences on nonuniform grids, we first derive an analytic approximation 
of a two-dimensional steady 
state vortex that nullifies the right-hand-side of~\eqref{eq:GPE}.
Then we perform a systematic comparison between the two approaches by measuring
the deviation of the numerical solution from the initial condition (being
steady the initial condition should remain preserved).
Finally, we explore the possibility to evaluate the solution obtained by
time-splitting Fourier methods on nonuniform grids designed to guarantee
higher spatial resolution in the proximity of 
vortex singularities.
\section{Accurate Pad\'e approximation of a 2d vortex}\label{sec:pade}
We seek a two-dimensional, steady-state solution of equation~\eqref{eq:GPE} 
that represents a straight vortex centered at the origin.
It is well-known that the classical two-dimensional Euler vortex of
circulation $\Gamma$ has azimuthal velocity
$u_\theta = \Gamma/(2\pi r)$ where $r=\sqrt{x_1^2+x_2^2}$ is the radius
and $\theta=\mathrm{atan2}(x_2,x_1)=\arg(x_1+\rmi x_2)\in(-\pi,\pi]$ is the azimuthal angle.
The Cartesian components of the velocity are thus
$u_1=-u_\theta \sin\theta = -\Gamma x_2/(2 \pi r^2)$ and 
$u_2=u_\theta \cos\theta =\Gamma x_1 /r^2$.
Therefore $\Vmy{u}=(u_1,u_2)=(\Gamma/(2 \pi))\Grad{\theta}$. 
This shows that the velocity field is solenoidal ($\Div{\Vmy{u}}=0$), that the
quantum mechanical phase, $S$, is simply the azimuthal angle $\theta$, and
that the quantum of circulation, in our dimensionless units, is
equal to $2 \pi$.  
In steady conditions, the continuity equation ensures that 
$\Div{(\rho\Vmy{u})}=0$, hence
$\Vmy{u}\cdot\Grad{\rho}=0$, which means that
$ \Grad{\rho}\cdot \Grad{\theta} = 0$. 
The solution $\rho = \bar\rho=\mathrm{const}$ leads to $\psi=\sqrt{\bar \rho}\rme^{\rmi \theta}$, which has infinite energy~\eqref{eq:energy} and must be 
rejected. 
The other possibility is that $\Grad{\rho} \perp \Grad{\theta}$. 
Since 
$\Grad\theta=\hat\theta/r$,
then $\Grad \rho$ is parallel to $\hat r$ and thus $\rho=\rho(r)$,
$\hat r$ and $\hat \theta$ being the unitary vectors in two-dimensional
polar coordinates.

In a two-dimensional domain we set
$\psi(x_1,x_2)=\rho(\sqrt{x_1^2+x_2^2})^{1/2}\rme^{\rmi \theta(x_1,x_2)}=f(\sqrt{x_1^2+x_2^2})\rme^{\rmi \theta(x_1,x_2)}$, where
$f(\sqrt{x_1^2+x_2^2})=f(r)$ is a function to be determined.
By imposing that $\psi$ is the steady solution of equation~\eqref{eq:GPE}, 
we find that $f$ satisfies the equation
\begin{equation}
f'' + \Fra{f'}{r} + f \left(1 - f^2 -\Fra{1}{r^2} \right)=0,
\label{eq:GPEf}
\end{equation}
with boundary conditions $f(0) = 0$, $f(\infty)=1$. 

Equation~\eqref{eq:GPEf} could be integrated numerically as it is, by
artificially bounding the infinite domain.
To avoid this problem, we resort to the change of 
variables $s=r/(1+r)$, $g(s)=f(r)$, which yields the equation
for $g(s)$
\begin{equation}\label{eq:g}
(s-1)^4g''+2(s-1)^3g'-\frac{(s-1)^3}{s}g'-\frac{(s-1)^2}{s^2}g+(1-g^2)g=0,
\end{equation}
defined in the finite domain $s\in(0,1]$, with boundary conditions
$g(0)=0$ and $g(1)=1$.
We solve equation~\eqref{eq:g} by central
second order finite differences with equally spaced discretization
points $s_i=i/N$, $i=1,2,\ldots,N$. 
Given the numerical solution $\tilde g$ of~\eqref{eq:g}, the numerical
approximation of the density is
\begin{equation}\label{eq:rhonum}
\rho_{\mathrm{num}}(r_i)=\left[\tilde g\left(\frac{r_i}{1+r_i}\right)\right]^2
\end{equation}
where $r_i=s_i/(1-s_i)$, $i=1,2,\ldots, N-1$. 
This rescaling provides denser points $r_i$ in the neighborhood of
the origin, where they are mostly needed 
(more than 95\% of the points $r_i$ are in the interval
$0\le r\le 20$). Nevertheless, the computation
of the initial solution for~\eqref{eq:GPE} on a two-dimensional grid, for
instance, requires $\rho_{\mathrm{num}}$ to be interpolated.

It would be therefore useful to have an analytic
approximation of $f(r)$.
However, since it is more convenient~\cite{B04} to find a Pad\'e approximation 
directly for $\rho(r)=\left[f(r) \right]^2$ rather than for $f(r)$, we rewrite 
equation~\eqref{eq:GPEf} in terms of $\rho(r)$ as
\begin{equation}
\rho''+\frac{\rho'}{r}-\frac{(\rho')^2}{2\rho}
-\frac{2\rho}{r^2}+2(1-\rho)\rho=0,
\label{eq:GPErho}
\end{equation}
with boundary conditions $\rho(0) = 0$, $\rho(\infty)=1$. 

It is known~\cite{B04,NW03} that Pad\'e approximations of $\rho(r)$ retain only
even degrees at both the numerator and denominator, that is
\begin{equation}
\rho(r) \approx
\frac{\sum_{j=0}^{p}a_j r^{2j}}{1+\sum_{k=1}^{q}b_k r^{2k}}=
\frac{a_0+a_1r^2+a_2r^4+\cdots+a_pr^{2p}}{1+b_1 r^2+b_2r^4+\cdots+b_qr^{2q}}.
\end{equation}
In order for this approximation to satisfy the boundary conditions, it must be
\begin{equation*}
\begin{aligned}
\rho(0)=0 & \implies  a_0 = 0\\
\rho(\infty)=1 & \implies  p = q, \quad b_q = a_p.\\
\end{aligned}
\end{equation*}
Given these simplifications, the diagonal Pad\'e approximation, with $2q-1$ coefficients
and both numerator and denominator of degree $r^{2q}$, is
\begin{equation}
\rho_q(r)=
\frac{a_1r^2+a_2r^4+a_3r^6+\cdots+a_qr^{2q}}{1+b_1 r^2+b_2r^4+b_3r^6+
\cdots+a_qr^{2q}}.
\end{equation}
In literature this approximation is normally limited to $q=2$~\cite{B04}, 
that is
\begin{equation*}
\rho_2(r)=\frac{a_1 r^2 +a_2 r^4}{1+b_1r^2+a_2r^4}
\end{equation*}
with 
\begin{equation*}
a_1=\frac{11}{32},\quad b_1=\frac{5-32a_1}{48-192a_1},\quad 
a_2=a_1\left(b_1-\frac{1}{4}\right).
\end{equation*}
Despite its widespread usage, this approximation is 
\emph{qualitatively wrong} (see, e.g.,~\cite{MRFL12}), 
because it reaches an unphysical maximum, above~$\rho(\infty) =1$, at 
$r_0=2\sqrt{6 \left( 4+3\sqrt{2} \right)}\approx14.065$, 
unique positive solution of $r^4-192r^2-1152$ obtained by imposing
$\rho'_2(r)=0$. 
The physical solution of equation~\eqref{eq:GPErho} must reach the limit value 
$\rho(\infty)=1$ monotonically, without overshooting.

Due to these limitations, we seek higher-order ($q>2$), 
monotonically increasing, Pad\'e expansions, namely
\begin{equation*}
\rho_3(r)=\frac{a_1r^2+a_2r^4+a_3r^6}{1+b_1r^2+b_2r^4+a_3r^6}
\quad\text{and}\quad
\rho_4(r)=\frac{a_1r^2+a_2r^4+a_3r^6+a_4r^8}{1+b_1r^2+b_2r^4+b_3r^6+a_4r^8}.
\end{equation*}

In order to determine the coefficients of a certain approximation $\rho_q(r)$, 
we compute the analytic expressions $\rho_q(r)$, $\rho'_q(r)$ and 
$\rho''_q(r)$ and substitute them in equation~\eqref{eq:GPErho} obtaining the 
form
\begin{equation}
\rho_q''+\frac{\rho_q'}{r}-\frac{(\rho_q')^2}{2\rho_q}
-\frac{2\rho_q}{r^2}+2(1-\rho_q)\rho_q=0 \iff
\frac{N_q(r)}{D_q(r)}=0.
\label{eq:GPErhoND}
\end{equation}
The numerator $N_q(r)$ is made of terms $r^{2k}$, which are in a number much
larger than the $2q-1$ degrees of freedom of the Pad\'e expansion.
For this reason equation~\eqref{eq:GPErho} cannot be satisfied exactly.
However, we can nullify the coefficients of $2q-1$ terms $r^{2k}$.
We can choose to start from higher- or lower-order coefficients in $N_q(r)$. 
We prefer to operate on lower-order powers of $r^{2k}$, i.e. $k=1, \dots, 2q-1$,
because we need a good approximation of $\rho_q(r)$ in a neighborhood of the
origin. 
Interestingly enough, we observe \emph{a posteriori} that canceling the
lower-order coefficients of $r^{2k}$ results in very small values of the 
coefficients of larger powers of $r$.
The step-by-step derivation of $\rho_q(r)$ for $q=2,3,4$ is reported 
in~\ref{sec:appendix}, whereas
tables~\ref{t:coeff23} and~\ref{t:coeff4} summarize all coefficients for the
expansions $\rho_q, q=2,3,4$.

\begin{figure}[!ht]
\centering
\includegraphics[scale=0.6]{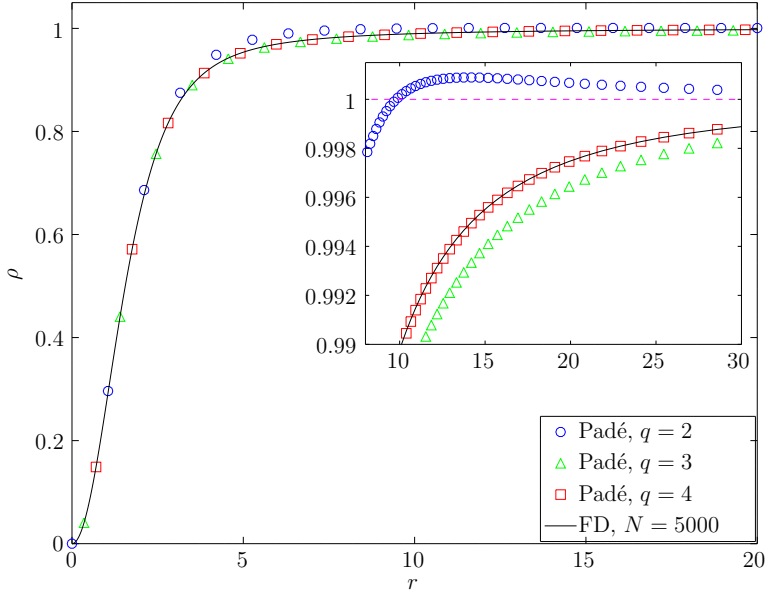}
\caption{Comparison between 
$\left[\tilde  g\left(\frac{r}{1+r}\right)\right]^2$, 
numerical solution of~\eqref{eq:g} obtained by second-order finite differences
on 5000 equispaced points,
and different Pad\'e approximations $\rho_q(r)$.}
\label{fig:rhopade}
\end{figure}

In Figure~\ref{fig:rhopade} we show the comparison between the 
numerical solution of~\eqref{eq:g} by employing
second-order central finite differences with 5000~points 
and different Pad\'e approximations $\rho_q(r)$ for $q=2,3,4$.
Visual inspection confirms that $q=2$ is a poor representation of the solution
of equation~\eqref{eq:GPErho}, especially for $4<r<20$.

\begin{figure}[!ht]
\centering
\includegraphics[scale=0.6]{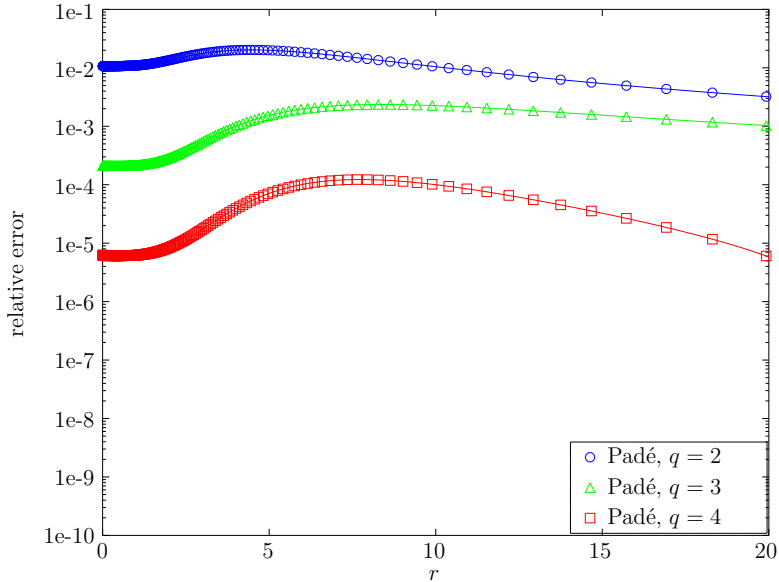}
\caption{Relative error between different Pad\'e approximations 
$\rho_q(r)$ and the numerical solution of~\eqref{eq:g} with
(equally spaced grid, second-order finite differences, 5000~points).}
\label{fig:errrho}
\end{figure}

To appreciate quantitatively the error with respect to the numerical solution,
in Figure~\ref{fig:errrho} we report the relative error in a semilog plot.
Interestingly, the relative error does not reach its maximum close to the 
origin,
meaning that any Pad\'e approximation reproduces quite well the behavior of the
vortex for $r \to 0$.
On the other hand, the maximum relative error is always reached for $r<10$, i.e.
in a region of interest for the numerical simulations that we will perform.
\section{Time splitting methods}\label{sec:splitting}
Widely used schemes for the numerical simulation of the dynamics 
of~\eqref{eq:GPE} are the so-called time-splitting methods and the finite
difference time domain methods (see~\cite{BC13} for a review).
If we restrict the options to second-order accurate schemes in time,
Time Splitting pseudoSPectral (TSSP) methods, Time Splitting Finite 
Difference (TSFD) methods
and Crank--Nicolson Finite Difference (CNFD) method conserve the mass at the
discretized level. 
However, CNFD is implicit and requires the solution
of a coupled nonlinear system at each time step.
For this reason we resorted to time splitting methods.
We refer the reader to~\cite{TCN09} for higher-order time
splitting methods.

In~\cite[Example 4.1]{BC13}
TSSP is suggested when the solution is smooth and TSFD otherwise, although
the hint comes from a one-dimensional numerical experiment. 
In what follows, we analyze two approaches: a classical time splitting Fourier 
method and a time splitting \emph{nonuniform} finite difference method.
In any case, equation~\eqref{eq:GPE} is split into two parts
\begin{subequations}
\begin{align}
\frac{\partial u}{\partial t}(t,x)&=\frac{\rmi}{2}\Lapl u(t,x)\label{eq:GPEkin}\\
\frac{\partial v}{\partial t}(t,x)&=\frac{\rmi}{2}\left(1-\abs{v(t,x)}^2\right)v(t,x)\label{eq:GPEpot}
\end{align}
\end{subequations}
where $x=(x_1,\ldots,x_d)\in\RR^d$.
The solution of the first equation depends on the space chosen for the 
discretization and will be described in the next two sections.
The second equation can be solved exactly, taking into account that $\abs{v}$
is preserved by the equation. Therefore
\begin{equation}
v(\tau,x)=\exp\left(\frac{\tau\rmi}{2}\left(1-\abs{v(0,x)}^2\right)\right)v(0,x)
\end{equation}
for any $x$ in the spatial domain. If we denote by 
$\rme^{\tau\mathcal{A}}u_n(x)$ and 
$\rme^{\tau\mathcal{B}(v_n(x))}v_n(x)$ the two partial numerical solutions, 
the approximation $\psi_{n+1}(x)$ of $\psi(t_{n+1},x)$, where $t_{n+1}=(n+1)\tau$, 
can be recovered by the so-called Strang splitting
\begin{equation*}
\begin{aligned}
\psi_{n+1/2}(x)&=\rme^{\tau\mathcal{A}}
\rme^{\frac{\tau}{2}\mathcal{B}(\psi_n(x))}\psi_n(x)\\
\psi_{n+1}(x)&=\rme^{\frac{\tau}{2}\mathcal{B}(\psi_{n+1/2}(x))}\psi_{n+1/2}(x).
\end{aligned}
\end{equation*}

\subsection{Time splitting Fourier method}
Equation~\eqref{eq:GPEkin} can be solved exactly in time within the Fourier 
spectral space.
A part from the error at machine-precision level coming from the
necessary direct and inverse Fast Fourier Transforms (FFTs),
the only possible considerable error might arise from an insufficient number 
of Fourier modes. This is usually not a big deal when 
approximating smooth solutions fastly decaying to zero, 
since spectral order of convergence takes place. 
For this to happen, the unbounded domain has to be truncated to a 
computational bounded domain $\Omega$ large enough to support the most 
of a periodic approximation of the solution. 
However, when simulating the dynamics of vortex solutions not decaying to zero,
as in our case where $\lim_{\abs{x}\to\infty}\abs{\psi(t,x)}=1$, there
are some issues to take into account: the low regularity of the
solution at the origin, due both to $\rho(t,0)$ and to $\theta(t,0)$,
and lack of periodicity at the boundaries, also considering the
usual extension of the computational domain and reflection
of the solution (see~\cite{KL93}). In fact, after 
such a mirroring, the solution takes the same values at opposite 
boundaries, but its derivatives do not.

In order to investigate the accuracy of Fourier approximation for vortex
solutions, we consider the Fourier series expansion of the function obtained
by mirroring

\begin{equation}\label{eq:psi0regular}
\psi_0(r,\theta)=(1-(1-\rho_4^{q/2}(r))\rme^{-r^2/\ell^2})\rme^{c\rmi\theta\rme^{-r^2/\ell^2}},
\quad \psi_0\colon[-L,L)^2\to\CC
\end{equation}
with respect to the axis $x=L$ and $y=L$. The final computational domain
is therefore $\Omega=[-L,3L)^2$, with $L=20$, 
for which $\rho_4(L,L)\approx0.99875$.
The choice of the parameters $q=1$, $\ell\to\infty$ and
$c=1$ corresponds to a two-dimensional straight vortex as described above.
Different choices provide more regular functions or functions 
fastly decaying to zero, for which the derivatives 
at the boundaries are almost periodic.
We compute a reference approximation by an expansion
into a series with $2048^2$ Fourier modes and compare it
with expansions
ranging between $M=m^2=16^2$ and $M=m^2=512^2$ modes, in the $L^2$ norm.
\begin{figure}[!ht]
\centering
\includegraphics[scale=0.6]{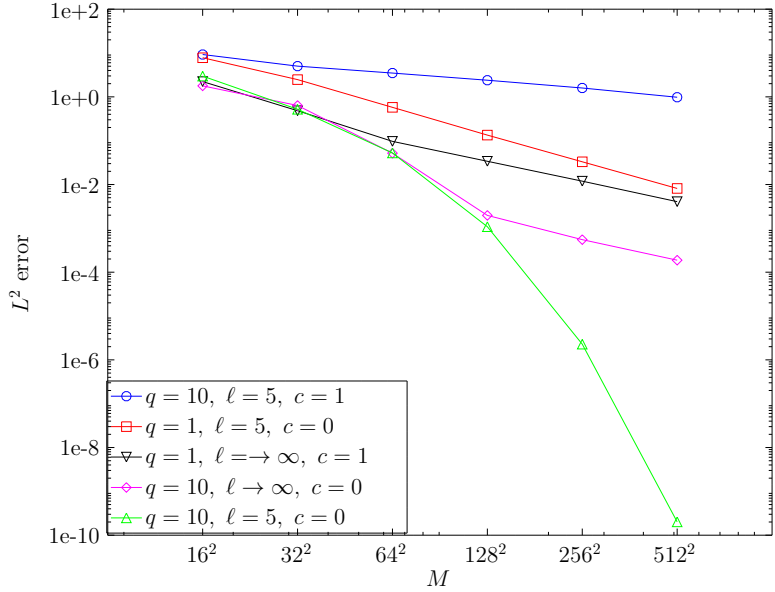}
\caption{Error behavior of the Fourier approximation of 
function~\eqref{eq:psi0regular} (extended to $[-L,3L)^2$ by mirroring) for different choices of the parameters.
Only the case of a regular and fast decaying to zero function 
($q=10$, $\ell=5$ and $c=0$, green, upward triangles) shows the typical
spectral rate of convergence.}
\label{fig:regularity}
\end{figure}
For a quite regular and periodic function, corresponding to $q=10$, $\ell=5$
and $c=0$ we observe in Figure~\ref{fig:regularity} a typical
spectral rate of convergence. For any other choice of the parameters, which
affects the regularity of the density ($q=1$, $\ell=5$, $c=0$), 
or the fast decay to zero of the function and its derivatives
($q=10$, $\ell\to\infty$, $c=0$) or the
regularity of the phase ($q=10$, $\ell=5$, $c=1$)
we observe a strong order reduction. The same reduction occurs for
the straight vortex ($q=1$, $\ell\to\infty$, $c=0$).

Increasing the number of Fourier coefficients so as to gain accuracy
is often not an option. In fact, due to the necessary mirroring, 
this corresponds to a huge growth of the degrees of freedom. 
Moreover, the use of hyperbolic sparse grids
(see~\cite{G07S}, for instance) is not possible, since the possibility
of discarding coefficients and grid points is given only for highly regular
solutions.

The low regularity of the solutions to be approximated and the needed 
duplication along axes in order to satisfy at least the periodicity of
the values of the solutions suggest to explore the alternative of a 
finite difference discretization in space.
\subsection{Time splitting finite difference method}\label{sec:TSFD}
The main advantage of a finite difference approach is that the mirroring of 
the solution is not required,
being the extension
of the bounded domain replaced by the imposition of homogeneous Neumann 
boundary conditions.

Given the low regularity of vortex solutions, we use centered second
order finite differences. 
With the aim of increasing the spatial resolution around
the vortex cores and keeping a reasonable degree of freedom,
we employ a set of nonuniform grid points (see~\cite{TA12}, for
instance, for locally adaptive finite element discretizations).

The discretization of the Laplace operator in one dimension 
with nonuniform finite
differences on $m$ points provides the nonsymmetric matrix
\begin{equation*}
A_1=
\begin{bmatrix}
-\frac{2}{h_1^2} & \frac{2}{h_1^2} & 0 & \ldots & 0\\
\frac{2}{h_1(h_1+h_2)} & -\frac{2}{h_1h_2} & \frac{2}{h_2(h_1+h_2)} & \ddots & 0\\
0 & \ddots & \ddots & \ddots & 0\\
0 & \ddots & \frac{2}{h_{m-2}(h_{m-2}+h_{m-1})} & -\frac{2}{h_{m-2}h_{m-1}} & \frac{2}{h_{m-1}(h_{m-2}+h_{m-1})}\\[1.8ex]
0 & \ldots & 0 & \frac{2}{h_{m-1}^2} & -\frac{2}{h_{m-1}^2}
\end{bmatrix}
\end{equation*}
where $h_i=x_{i+1}-x_i$, $x_1=-L$, $x_m=L$. 
This is not exactly a second order approximation, although
a discretization in which $h_{i+1}=(1+\delta)h_{i}$ and a refinement with
$h_{j+1}=(1+\delta)^{1/2}h_{j}$ yields a first order term in the error decaying
faster than the second order term (see~\cite[\S~3.3.4]{FP02}). The
approximation for the two-dimensional and the three-dimensional cases
can be simply obtained by Kronecker products with the identity matrix.
If we call $A$ the corresponding matrix, equation~\eqref{eq:GPEkin} is 
transformed into the system of ordinary differential equations 
\begin{equation}\label{eq:ODE}
y'(t)=\frac{\rmi}{2} Ay(t),\quad y(t)\in\CC^{M\times 1}.
\end{equation}

Given the importance of the mass preservation, we investigate this issue
for the numerical solution of system~\eqref{eq:ODE}.
A quadrature formula with positive
weights for the computation of
the mass writes
\begin{equation*}
\int_\Omega \abs{\psi(t,x)}^2\de x\approx w^T\abs{y(t)}^2,\quad w\in\RR^{M\times
1}_{+}.
\end{equation*}
It can be written as 
\begin{equation*}
y(t)^*Wy(t)
\end{equation*}
where $y(t)^*\in\CC^{1\times M}$ denotes the transposed conjugate vector of $y(t)$
and $W$ the matrix with diagonal $w$.
We define $z(t)=W^{1/2}y(t)$ such that
\begin{equation}\label{eq:ODEsh}
z'(t)=\frac{\rmi}{2}A_wz(t)
\end{equation}
with $A_w=W^{1/2}AW^{-1/2}$. 
If $A_w$ is \emph{symmetric}, then
the solution $z(\tau)=\exp(\tau\rmi/2 A_w)$ is an orthogonal matrix and 
\begin{equation*}
z(\tau)^*z(\tau)=z^*(0)z(0).
\end{equation*}
This means that
\begin{equation*}
\begin{split}
y(\tau)^*Wy(\tau)&=(W^{1/2}y(\tau))^*(W^{1/2}y(\tau))=z(\tau)^*z(\tau)=z(0)^*z(0)=\\
&=(W^{1/2}y(0))^*(W^{1/2}y(0))=y(0)^*Wy(0)
\end{split}
\end{equation*}
and therefore system~\eqref{eq:ODE} preserves the mass at the discrete level
if $W$ makes $A_w$ symmetric. 
From the structure of the matrix $A_1$, it is clear that the vector of trapezoidal weights
$w^T=[h_1,h_1+h_2,h_2+h_3,\ldots,h_{m-1}]$ gives a matrix $W_1$ such that 
$W_1A_1$ is symmetric. The extension to $W$ in the two-dimensional and
three-dimensional cases is
trivial 
and this is enough to get $A_w$ symmetric as well, in fact
$A^TW=WA\iff W^{-1/2}A^TW=W^{1/2}A\iff W^{-1/2}A^TW^{1/2}=A_w^T=W^{1/2}AW^{-1/2}=A_w$.
We conclude that equation~\eqref{eq:ODE} preserves the mass at the 
discrete level  whenever
the trapezoidal rule is used as quadrature formula and this is
easily extended to any space dimension.

System~\eqref{eq:ODEsh} could be solved, for instance, by the Crank--Nicolson
scheme
\begin{equation*}
z_{n+1}=z_n+\frac{k\rmi}{4}A_w z_n+\frac{k\rmi}{4}A_w z_{n+1},
\end{equation*}
which preserves the discrete mass
being $A_w$ symmetric (see~\cite{BC13}). 
This scheme is second order accurate in time, therefore the size
of the time step $k$ has to be chosen such that the
error is smaller than the time splitting error.
Moreover, Crank--Nicolson scheme requires the solution of a linear system
of equations with matrix $(I-k\rmi A_w/4)$ at each time step $k$.
Although this is not a big deal in one space dimension, since the
matrix is tridiagonal, in higher dimensions the discretization 
yields a large, sparse, complex symmetric matrix.
This implies the use of preconditioned 
Krylov solvers for general matrices such as
GMRES or BiCGStab or minimal residual methods for complex symmetric
systems (see~\cite{C13}). Iterative methods converge to the solution up to a 
specified tolerance which therefore influences the mass conservation and
the whole accuracy of the result.
Given these complications, we prefer to consider a direct approximation
of the exact solution 
\begin{equation*}
z_{n+1}=\exp(\tau\rmi/2 A_w)z_n.
\end{equation*}
Nowadays there are several options for the computation of the
action of the matrix exponential to a vector. We refer to~\cite{CKOR14}
for a review of \emph{polynomial} methods which do not require the
solution of linear systems. In this way, the kinetic linear part~\eqref{eq:GPEkin}
is solved exactly in time, as in the Fourier spectral method.
\section{Numerical experiments}
In Section~\ref{sec:pade} we have derived various approximations of $\rho(r)$ 
for a straight, two-dimensional vortex, whose 
wavefunction is $\psi(r,\theta)=\sqrt{\rho}\rme^{\rmi\theta}$.
In order to quantitatively compare the two methods introduced in
Section~\ref{sec:splitting}, we measure the 
preservation of such a steady solution by reporting the relative error 
\begin{equation}
\max_{0<\abs{r}\le R}\frac{\abs{\psi_n(r,\theta)-\psi_0(r,\theta)}}
{\abs{\psi_0(r,\theta)}},
\quad n=1,2,\ldots,T/\tau
\label{eq:error}
\end{equation}
with
$\psi_0(r,\theta)=\abs{\psi_0(r)}\rme^{\rmi\theta}$,
where $\abs{\psi_0(r)}$ is either $\sqrt{\rho_q(r)}$ or 
$\sqrt{\rho_{\mathrm{num}}(r)}$, the latter evaluated at any required $r$
by linear interpolation
of~\eqref{eq:rhonum}.
The origin is excluded since $\psi_0$ is zero therein.
The time step $\tau$ is chosen such that $T/\tau$ is an integer, where
$T$ is the final simulation time. 
In all our experiments, we selected $T=10$, a reasonable value in quantum 
fluids simulations~\cite{ZCBB12,AZCPPB14,ZR15}. 
The maximum over the continuum set $\{0<\abs{r}\le R\}$ in the error above 
is approximated by the maximum over a discrete set which will be specified
later. 

Although the preservation of the initial state may seem a trivial test, 
it is in fact a reliable and necessary experiment in order to validate 
the effectiveness of the proposed numerical methods. Thanks to the reliability of
the analytic solution, this test can show the influence of both the spatial 
approximation and the time splitting error in the numerical discretization 
of the PDE~\eqref{eq:GPE}.

In what follows we will employ either TSSP (Fourier) or TSFD.
For a computational grid with $m\times m$ grid points 
in the physical domain of interest, TSSP requires a total of $M=2 m \times 2m 
= 4 m^2$ degrees of freedom due to mirroring, whereas TSFD requires only 
$M=m^2$ degrees of freedom thanks to homogeneous Neumann boundary conditions.

\subsection{Comparison between different approximations of the initial condition}
\begin{figure}[!ht]
\centering
\includegraphics[scale=0.3,bb=0 360 565 795]{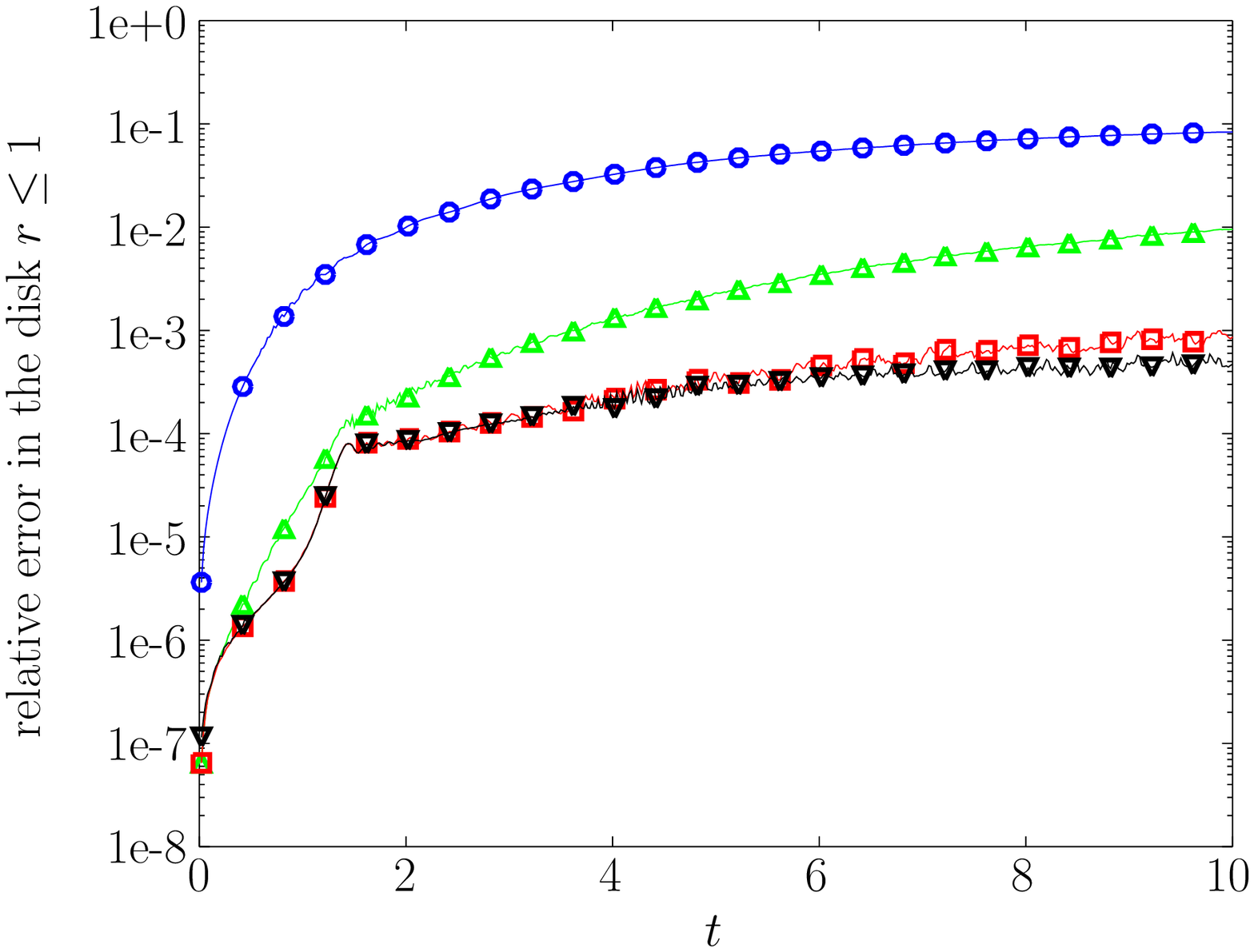}
\includegraphics[scale=0.3,bb=0 360 565 795]{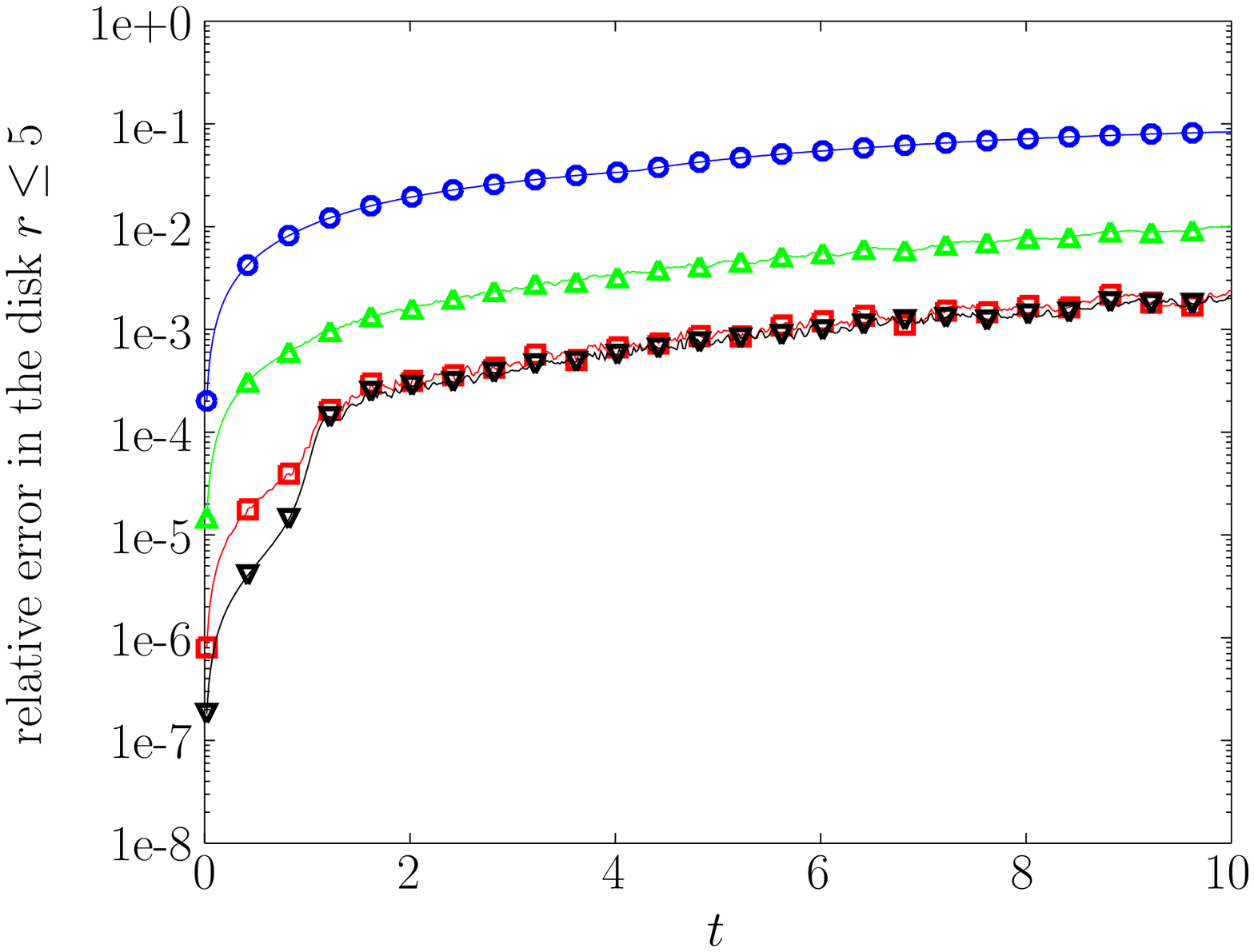}\\
\includegraphics[scale=0.3,bb=0 360 565 795]{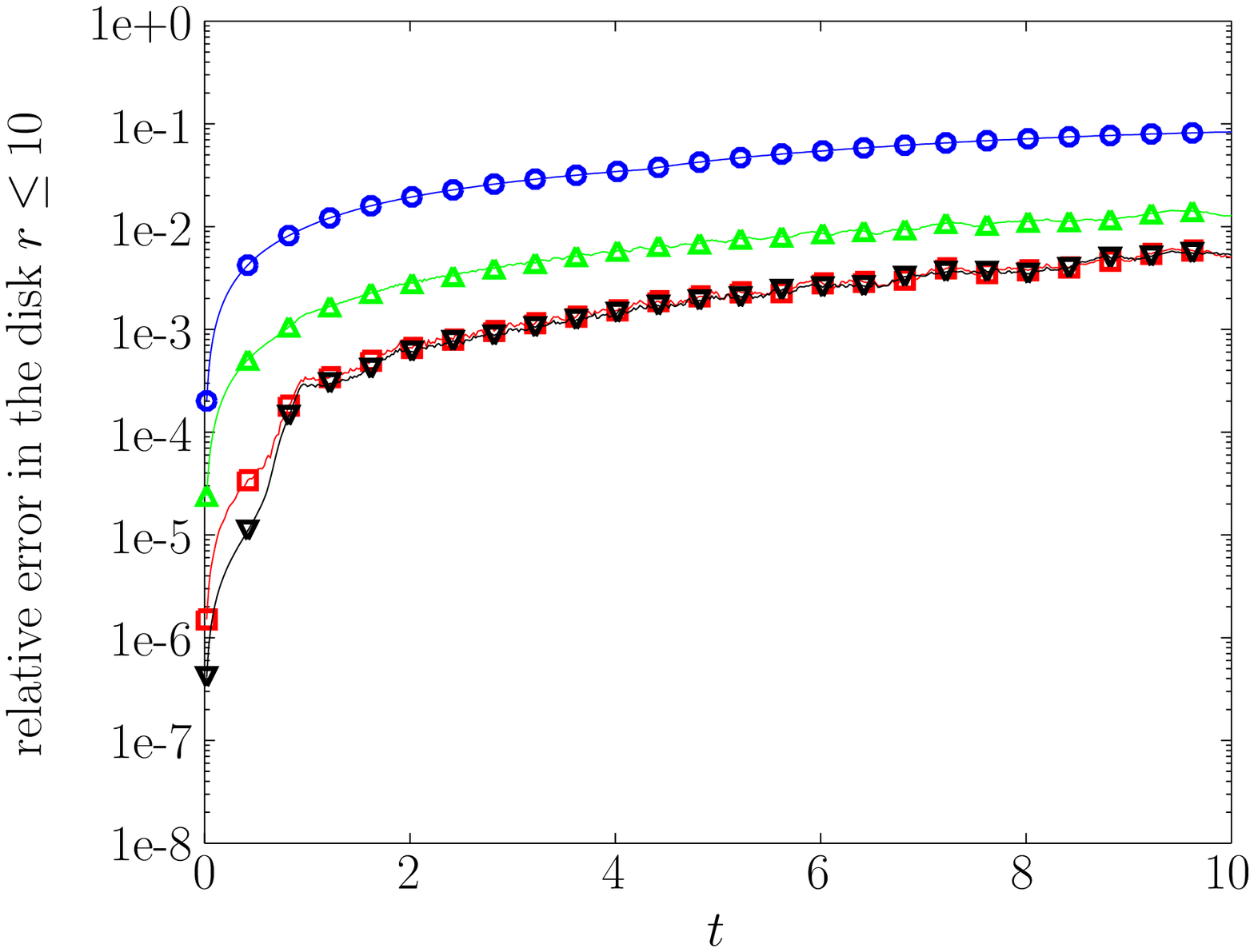}
\includegraphics[scale=0.3,bb=0 360 565 795]{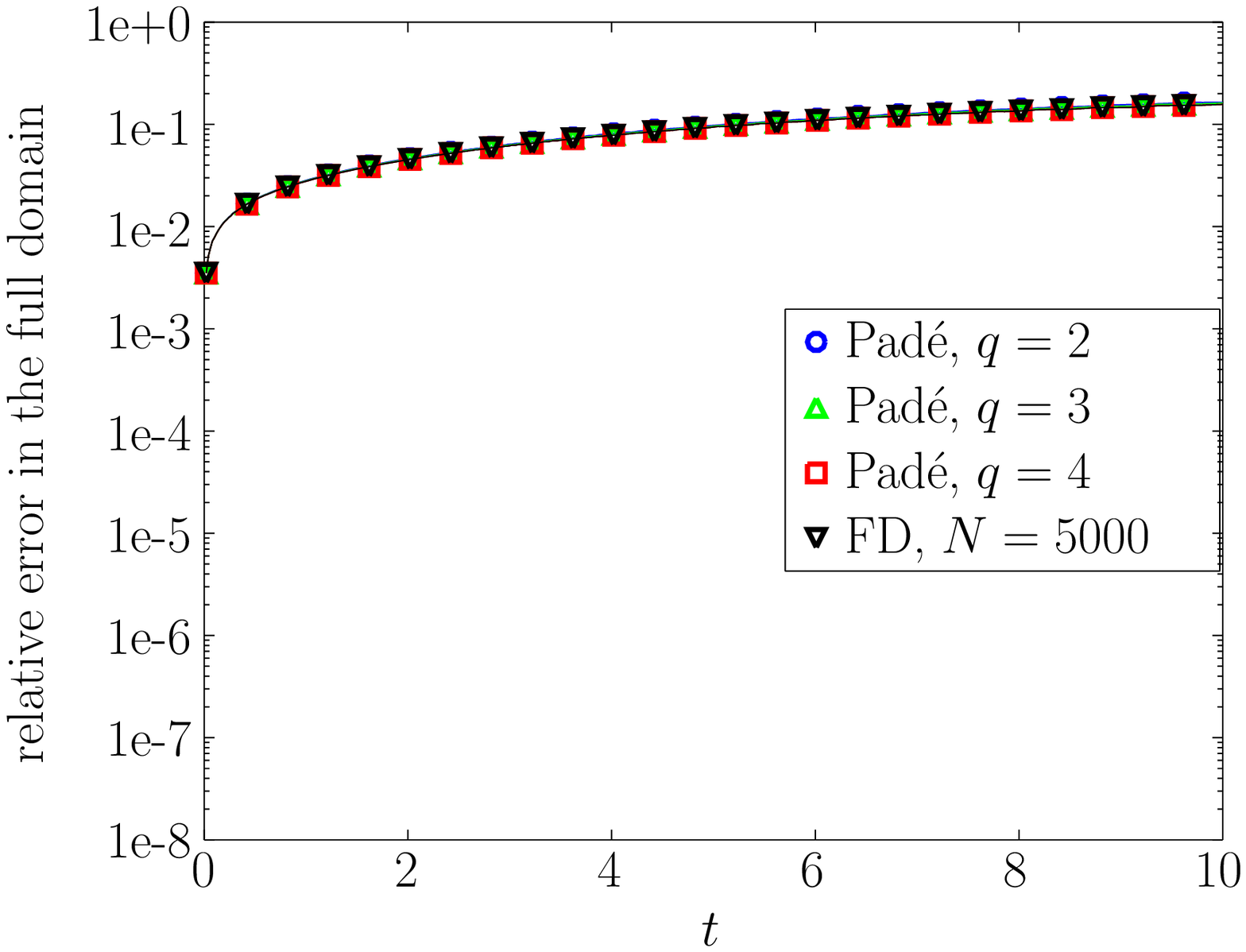}
\caption{Comparison of the relative error as defined by~\eqref{eq:error} for
different choices of the initial condition, Fourier approach.}
\label{fig:diffIC}
\end{figure}
We preliminary test the reliability of the three Pad\'e approximations 
$\rho_q, q=2,3,4$ 
and the numerical solution of equation~\eqref{eq:g} obtained 
by central second order finite differences with $N=5000$ uniformly distributed 
discretization points.
For the solution of the GPE~\eqref{eq:GPE} 
we employ TSSP with Fourier basis functions on a uniform two-dimensional
computational grid.
For this reason, the numerical solution $\tilde g(s)$ of equation~\eqref{eq:g}
must be interpolated.

Results are reported in Figure~\ref{fig:diffIC}, where the relative error
defined by~\eqref{eq:error} is plotted versus time for different disks.
The number of Fourier modes is fixed to $m=2\cdot 200$, i.e. 
$M=1.6 \times 10^5$ degrees of freedom.
We compare the solution at each time step with the initial condition 
on the grid nodes within the considered disk.
The worst approximation of the steady-state solution is the commonly used
$\rho_2$ Pad\'e approximation, whereas 
$\psi_0(r,\theta)=\sqrt{\rho_4(r)}\rme^{\rmi \theta}$ turns out to be
as accurate as the numerical solution.
For this reason,  in the following experiments we will 
consider only
$\psi_0(r,\theta)=\sqrt{\rho_4(r)}\rme^{\rmi \theta}$.
All curves collapse on each other in the case of the largest disk, meaning that 
the maximum error occurs at the boundaries, mainly due to the non-periodicity 
of the solution.

\subsection{Uniform vs. nonuniform finite differences}
\begin{figure}[!ht]
\centering
\includegraphics[scale=0.3,bb=0 360 565 795]{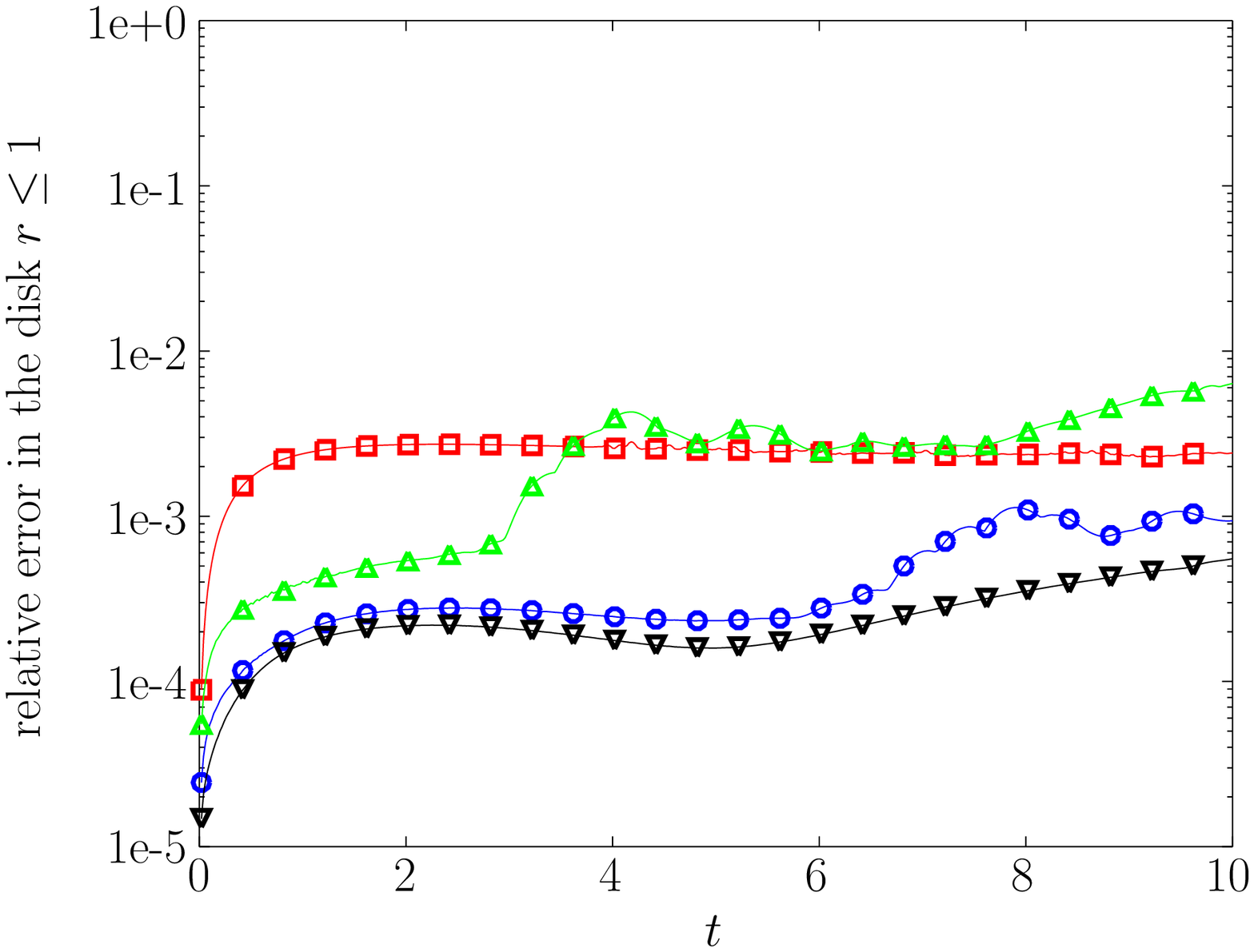}
\includegraphics[scale=0.3,bb=0 360 565 795]{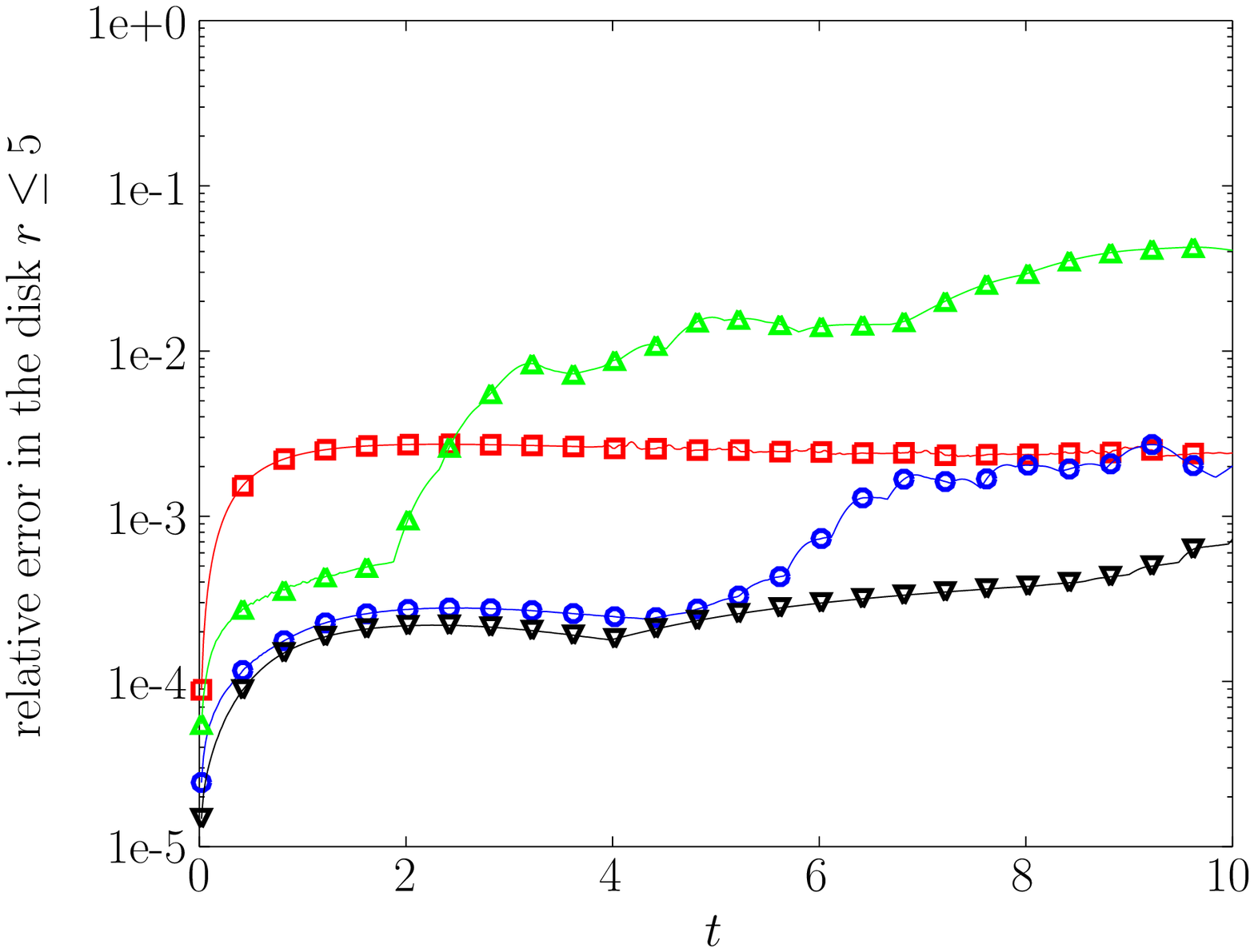}\\
\includegraphics[scale=0.3,bb=0 360 565 795]{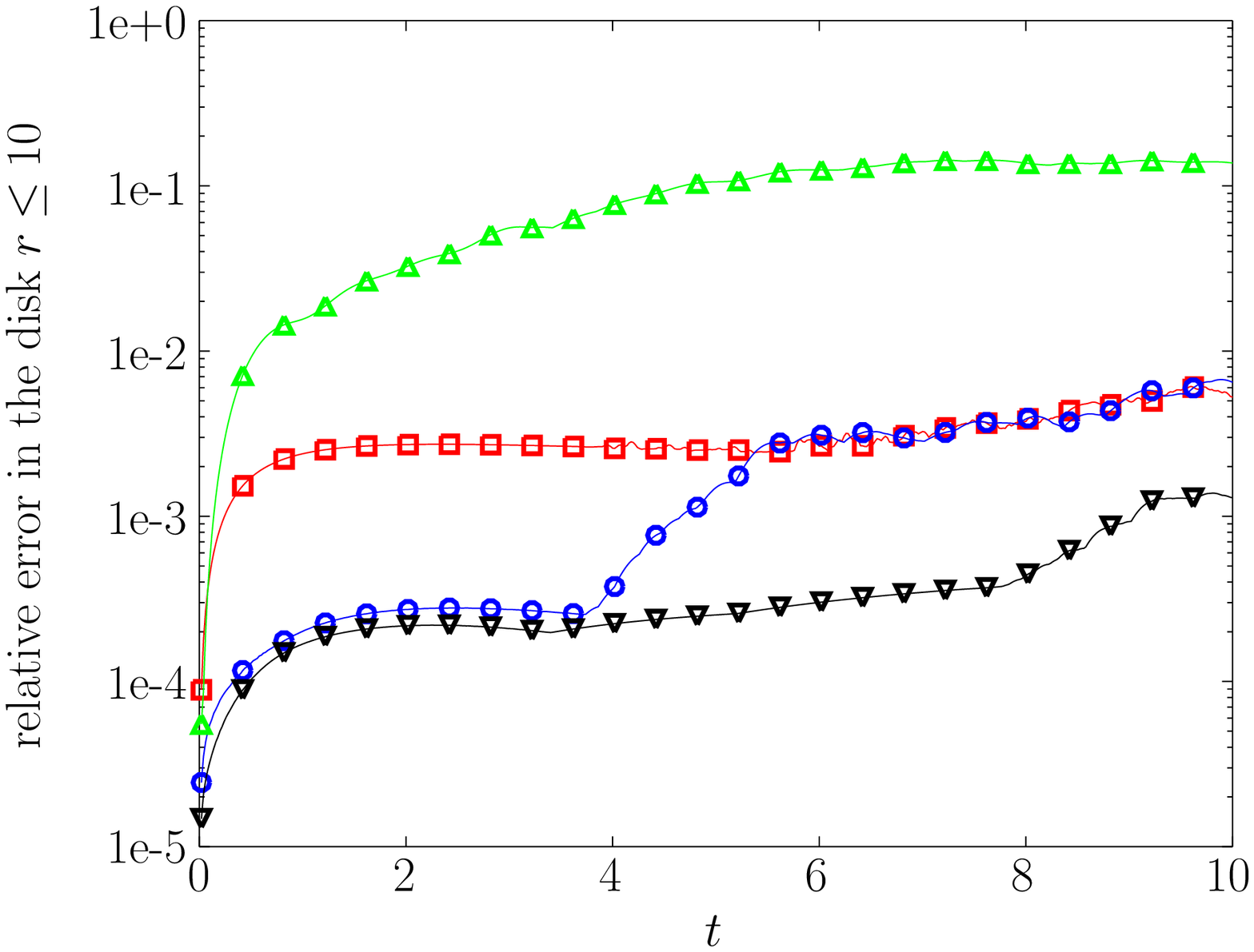}
\includegraphics[scale=0.3,bb=0 360 565 795]{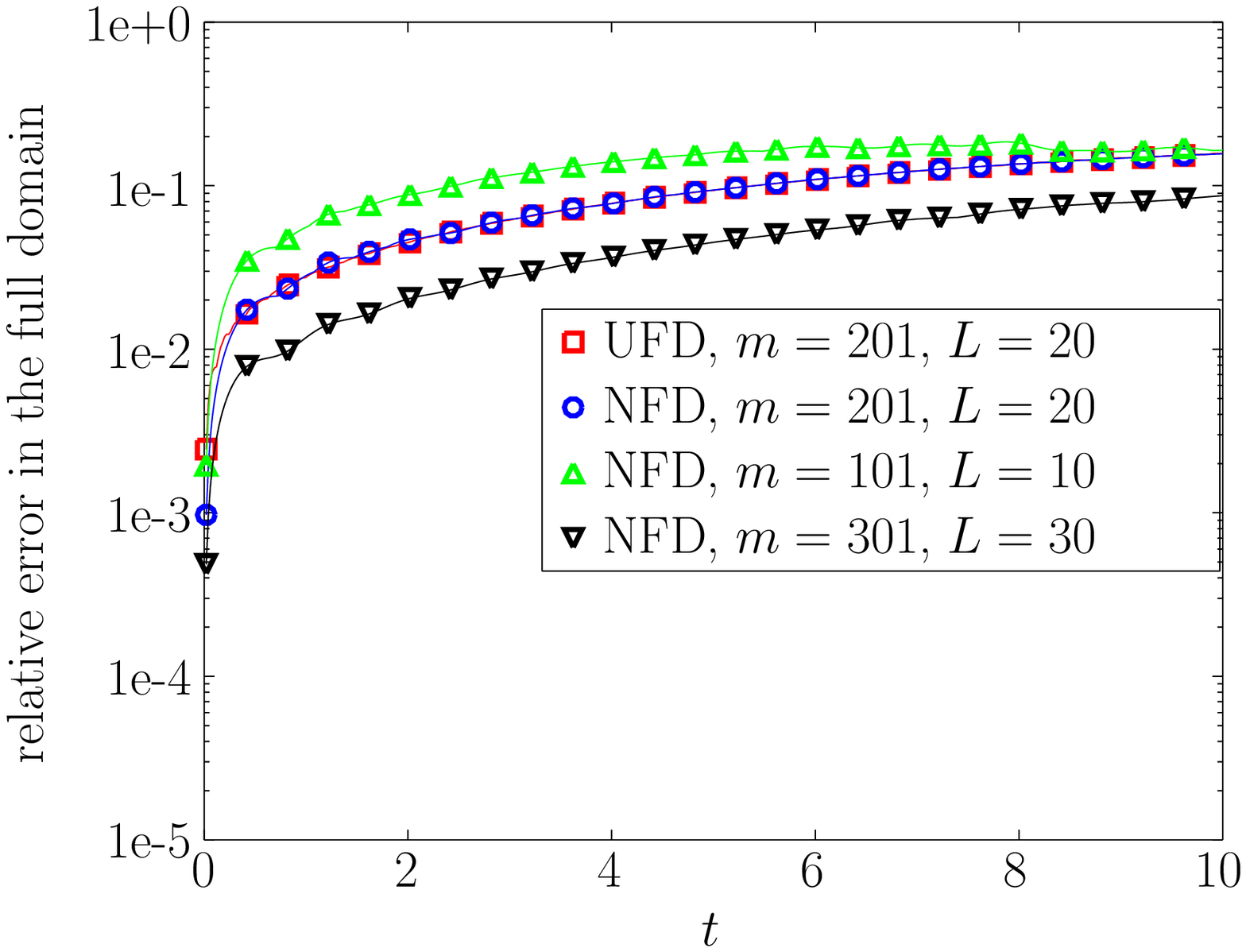}
\caption{Comparison of the relative error as defined by~\eqref{eq:error} for
central finite differences on uniform (UFD) and nonuniform (NFD) grids.}
\label{fig:univsnon}
\end{figure}
We compare the performance of the finite difference approximation on a uniform
versus nonuniform grid.
The uniform grid has the same step-size $h$ as the grid 
employed for the Fourier approach and reported in Figure~\ref{fig:diffIC}.
The nonuniform grid is generated by taking into account different constraints.
Given the smallest step-size $h_\mathrm{min}=h_1=0.05$ at the origin (in the vortex core), we
linearly increase the step-size according to 
$h_{i+1}=(1+\delta)h_{i}$ in both $x$ and $y$
and in both positive and negative directions.
We choose $\delta$ so as to reach the boundaries exactly, and in order to 
keep the ratio $K=h_\mathrm{max}/h_1 \approx 10$, where
$h_\mathrm{max} = \max h_i$. 
The number of points of the nonuniform grid is chosen such that
the mean value of $\{h_i\}$ equals the step-size of the uniform grid.
These constraints guarantee a reasonable nonuniform grid.

The comparison between the uniform and nonuniform grids is shown in 
Figure~\ref{fig:univsnon}.
In all cases the initial condition is
$\psi_0(r,\theta)=\sqrt{\rho_4(r)}\rme^{\rmi \theta}$, thus
the error for the uniform grid can be compared directly with that in
Figure~\ref{fig:diffIC} for the case of Pad\'e approximation with $q=4$ (red
squares in both Figures).

We first focus on the results with the same number of points and the same 
boundaries, i.e. red squares and blue circles in Figure~\ref{fig:univsnon}.
The discrete mass variation along time integration is comparable and of
order $10^{-13}$.
This confirms the conservation of mass also for the case of nonuniform grid,
as discussed in Section~\ref{sec:TSFD}.
In the nonuniform case the error is roughly one order of magnitude smaller than 
in the uniform case
on small disks and for $t$ not too large, whereas the 
curve of the nonuniform case tends to jump onto the uniform one after
a certain time as the radius of the disk of interest increases.
This suggests the idea that the error arises at the boundaries, where $\psi_0(r,\theta)$
does not exactly fulfill Neumann boundary conditions.

Motivated by this, we have changed the boundaries from $L=20$ to $L=30$ and
$L=10$ to check the dependency of the error on the choice of the truncated
domain.
In doing so, we have preserved the constraints on the nonuniform grids 
discussed above, obtaining $M=101^2$ degrees of freedom for $L=10$ 
(upward green triangles) and $M=301^2$ degrees of freedom for $L=30$ 
(downward black triangles).
With reference to Figure~\ref{fig:univsnon},
the domain bounded at $L=10$ is clearly too small and the error is always very
large compared to all the other cases.
On the other hand, the curves for $L=20$ (blue circles) and $L=30$ (downward
black triangles) behave roughly in the same way up to a certain 
value of $t$, after which the case $L=20$ consistently show larger errors than
the case $L=30$.
This reinforces the claim that the error arises from the borders.

\subsection{Comparison between Fourier spectral method and nonuniform 
finite differences}
\begin{figure}[!ht]
\centering
\includegraphics[scale=0.3,bb=0 360 565 795]{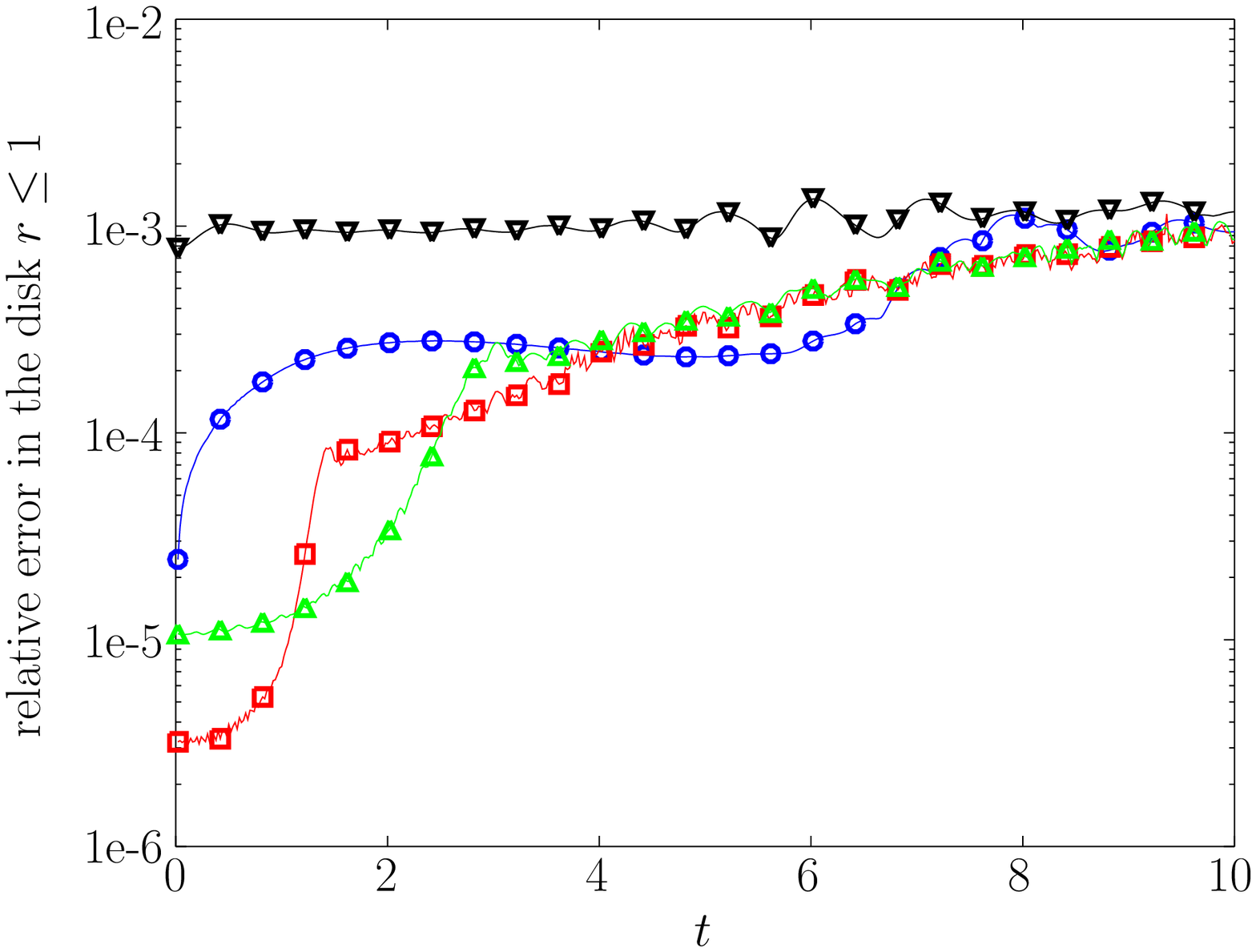}
\includegraphics[scale=0.3,bb=0 360 565 795]{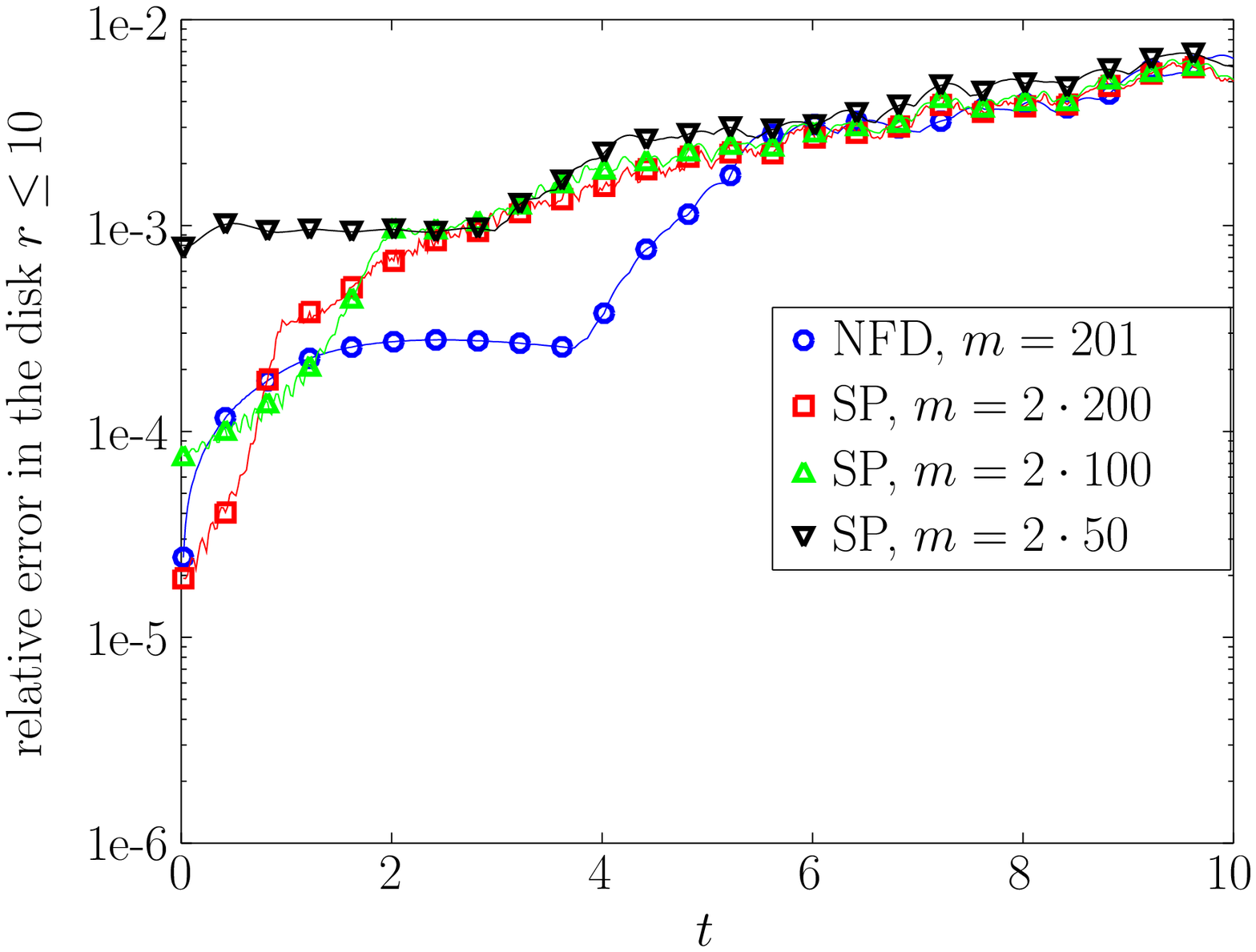}
\caption{Comparison of the relative error as defined by~\eqref{eq:error} between
nonuniform finite differences~(NFD) and spectral Fourier~(SP) for
different numbers of Fourier modes.}
\label{fig:TSFDvsTSSP}
\end{figure}

Now we concentrate on our main goal, which is the comparison between TSSP with
Fourier basis function on uniform grids and TSFD on a nonuniform grid that we
fix to $h_\mathrm{min}=0.05$, $L=20$, $m=201$.
In order to compare  the error defined by~\eqref{eq:error} for the two methods,
we always evaluate the TSSP solution on the nonuniform grid points (spectral
solutions can be evaluated everywhere).
This set of points has the advantage of being denser in the
vortex core, where higher spatial resolution is desirable.
Results are reported in Figure~\ref{fig:TSFDvsTSSP}, where SP stands for
spectral and NFD for nonuniform finite differences.

Keeping in mind that the spectral Fourier approach needs mirroring, i.e. the
number of modes in each direction must be doubled,
we first choose a number of Fourier modes $m=2\cdot 200$ in each direction 
to make it equal to the number of points of the 
reference case for nonuniform finite differences ($m=201$) in the
physical (un-mirrored) domain.
The overall behavior of the error for these two cases is comparable: 
TSSP (red squares) performs better than TSFD (blue circles) for small 
values of $t$,  whereas the opposite
happens for intermediate values of $t$.
For large $t$ the two curves collapse on each other.

Due to the fact that TSSP needs mirroring, i.e.
$M_\mathrm{TSSP}=4M_\mathrm{TSFD}$, in Figure~\ref{fig:TSFDvsTSSP} we explore 
also the cases with less Fourier modes, namely
$m=2\cdot 100$ (upward green triangles) and 
$m=2\cdot 50$ (downward black triangles).
As observed for the case $m=2\cdot 200$, in the long term all curves seem to
provide similar errors, regardless of the disk radius.
On the other hand, for small values of $t$, the number of degrees of freedom
plays a r\^ole in that a larger number of Fourier modes ensures smaller errors.

It is important to keep in mind that, for what seen in
Figure~\ref{fig:regularity}, the TSSP Fourier approach does not retain the
spectral accuracy because of the singular nature of the solution at the origin
and the lack of periodicity at the boundaries.

As a final remark, we observe that the error of the Fourier solution computed 
on its own uniform grid, reported in Figure~\ref{fig:diffIC} with red squares, 
is smaller than the error of the Fourier solution evaluated on the nonuniform 
grid, reported in Figure~\ref{fig:TSFDvsTSSP} with red squares.

\subsection{Maximum resolution of Fourier spectral method}
As expected, from Figure~\ref{fig:TSFDvsTSSP} we have seen that the smaller the number 
of Fourier modes, the larger the relative error with respect to the initial
condition.

\begin{figure}[!ht]
\centering
\includegraphics[scale=0.3,bb=0 360 565 795]{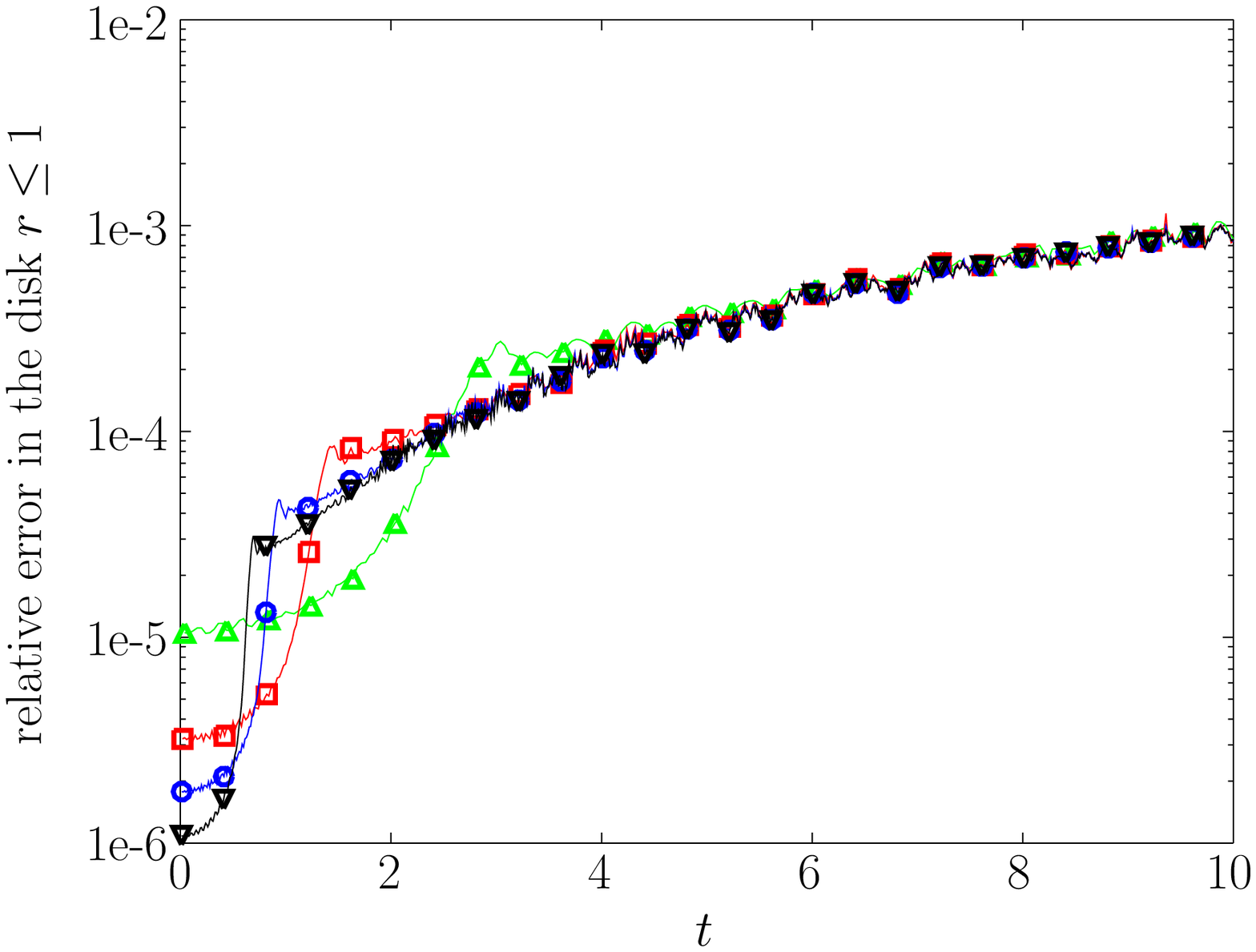}
\includegraphics[scale=0.3,bb=0 360 565 795]{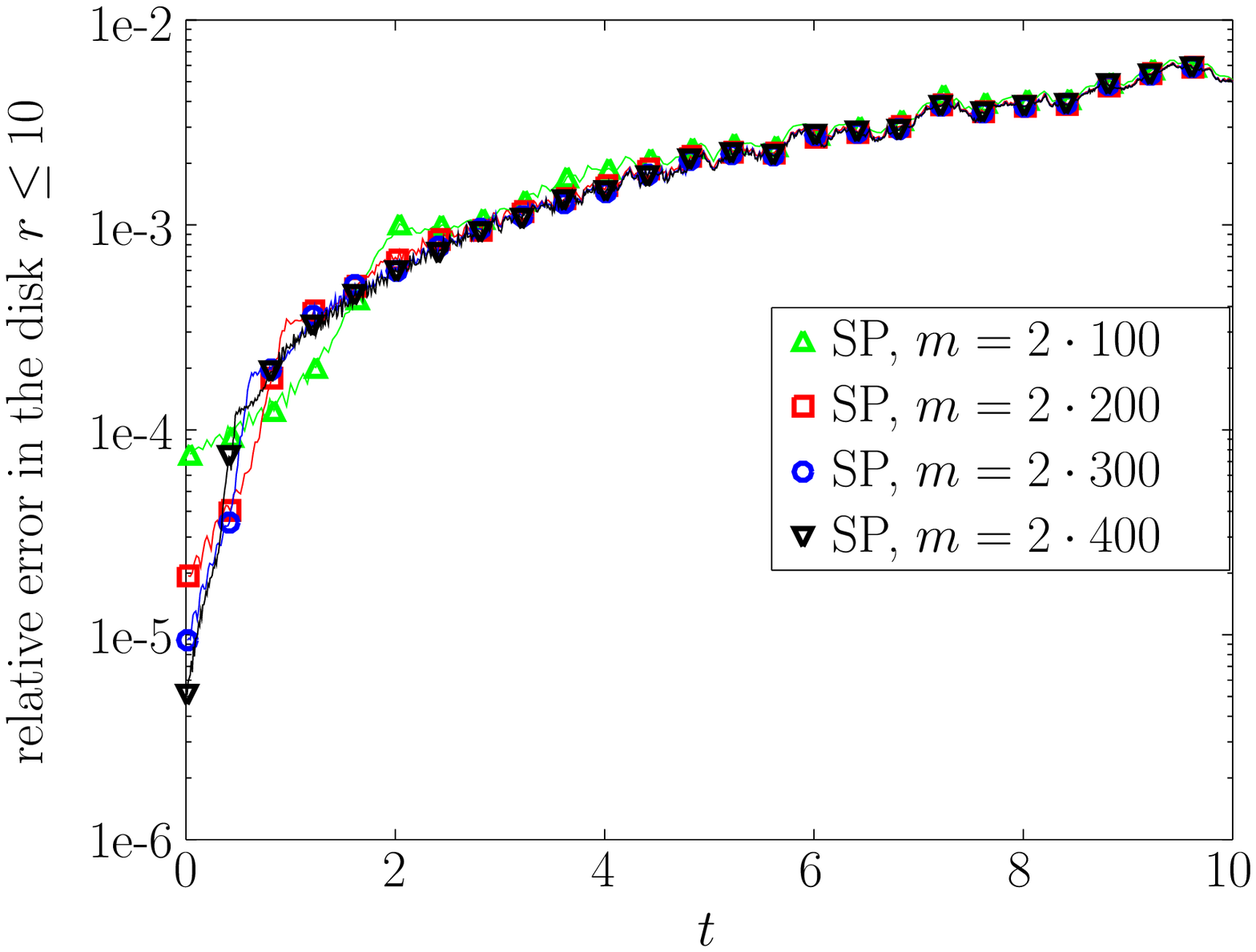}
\caption{Comparison of the relative error as defined by~\eqref{eq:error} 
for increasing number of Fourier modes.}
\label{fig:maxresTSSP}
\end{figure}

We wish to check if there exists an upper limit to the maximum resolution of 
Fourier spectral method.
For doing so, we increase the number of Fourier modes and, proportionally, the 
number of time steps as suggested in~\cite{BJM02}.
Results are shown in Figure~\ref{fig:maxresTSSP}.
We observe high accuracy in the core (see smaller disk, left 
plot) for small values of $t$, immediately followed by saturation.
In a larger disk (right plot), saturation kicks in almost immediately.
The errors reported in Figure~\ref{fig:maxresTSSP} suggest that $m=2\cdot 200$ 
is a reasonable value of Fourier modes for the preservation of a 
two-dimensional quantum vortex.

\subsection{Fourier evaluation on nonuniform grids}
As explained in the Introduction, our motivation to explore the nonuniform
finite difference approach is based on the need, for the study of vortex 
reconnections~\cite{ZR15}, of high \emph{local} spatial resolution and, 
possibly, accuracy.
\begin{figure}[!ht]
\centering
\includegraphics[scale=0.3,bb=0 360 565 795]{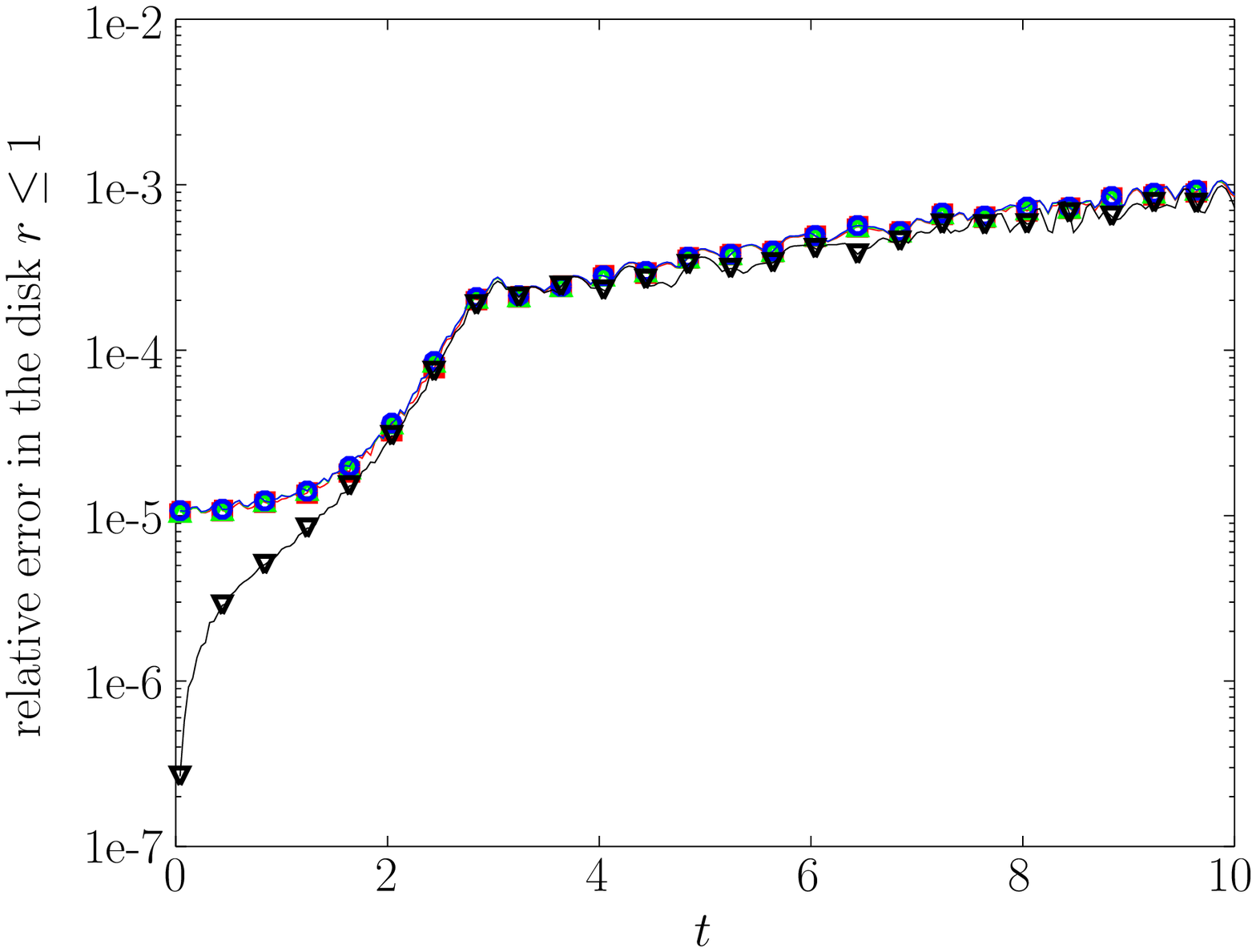}
\includegraphics[scale=0.3,bb=0 360 565 795]{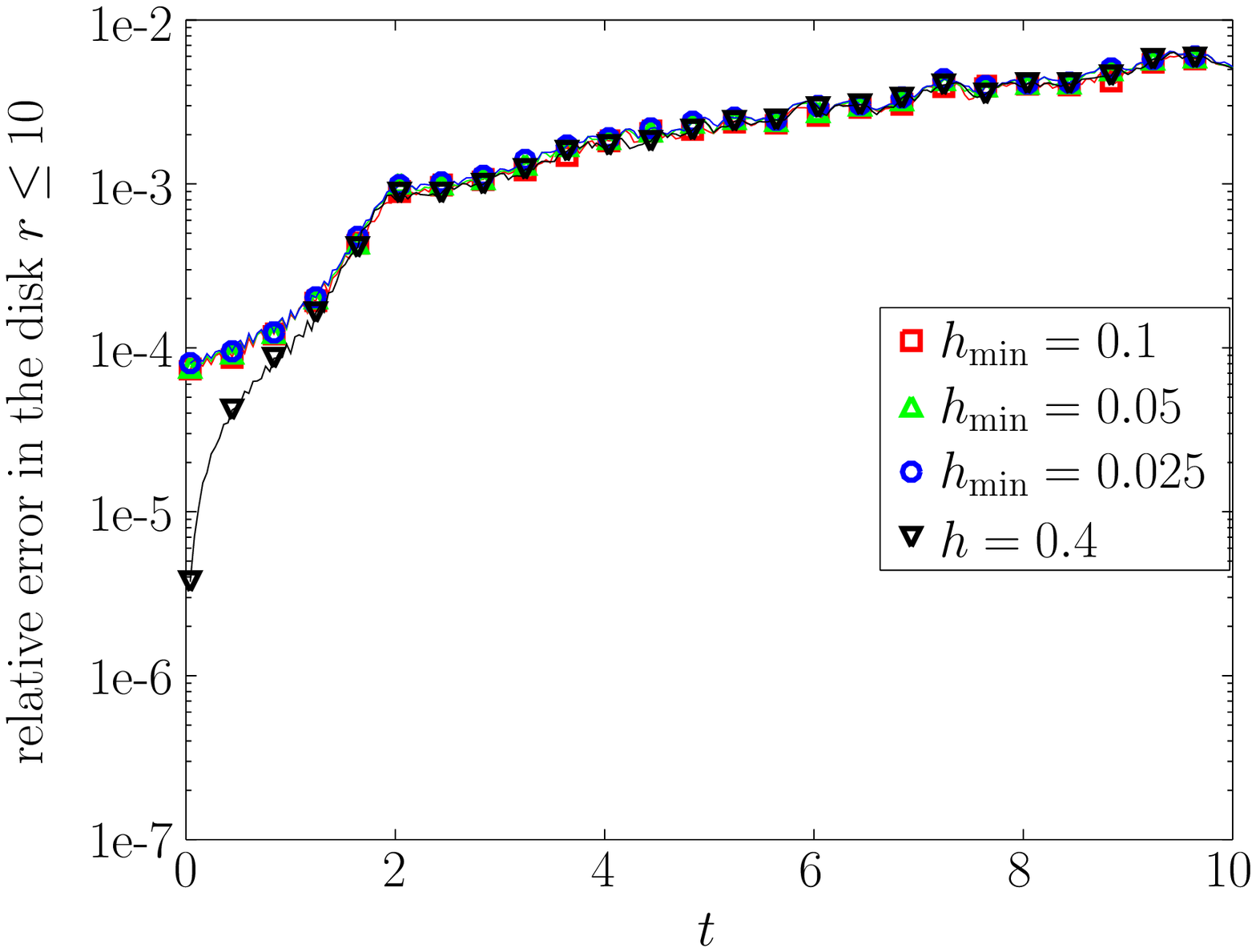}
\caption{Comparison of the relative error as defined by~\eqref{eq:error} 
for different evaluations of the Fourier solution ($m=2\cdot 100$)
at nonuniform grids.}
\label{fig:TSSPcoll}
\end{figure}

Instead of increasing the number of Fourier modes so as to reach higher 
\emph{global} spatial resolution, one can resort to a TSSP method with a 
reasonable number of modes (considering that mirroring is needed) and then 
evaluate the TSSP solution on a nonuniform grid, with denser points where they
are needed.
In Figure~\ref{fig:TSSPcoll} we compare the reasonable case $m=2\cdot 100$, for
which the number of modes is relatively small, but not too small, with different
nonuniform grids. 
We notice that the numerical integration itself is carried out only once 
and the Fourier coefficients of the solution are stored at each time step.
The evaluation at the grid points is performed afterward, in the post-processing 
stage, as many times as desired. Moreover,
tools like the Nonuniform Fast Fourier Transform (NFFT, see~\cite{KKP09}) can
be employed for the fast evaluation of trigonometric polynomials at arbitrary
point sets.
The constant spatial step-size of the TSSP method is $h=0.4$, whereas $h_\mathrm{min}$
stands for the minimum value of the step-size, in the proximity of the origin,
for the nonuniform grids.
As seen before, there is a substantial difference in the error only for $t<2$,
whereas  for larger values of $t$ evaluating the Fourier solution on a 
nonuniform grid does not worsen the solution.
It is important to note that vortex reconnections, usually, require a dynamics
that takes a time of at least $t=10$.
Evaluating a TSSP solution on a locally refined grid is, thus, a very promising
approach to study quantum vortex reconnections.
\section{Conclusions}
After deriving a new accurate Pad\'e approximation for the density distribution 
of a two-dimensional steady-state vortex, we have used it as the initial 
condition for the Gross--Pitaevskii equation to test the performance of the 
time-splitting Fourier method.
Although it cannot retain its classical spectral accuracy in space, 
being as accurate as low-order finite difference on nonuniform grids, 
it preserves quite well the steady-state solution, especially in the 
neighborhood of the singularity. 
The advantage of a post-processing evaluation on arbitrary points makes 
this approach suited for applications where \emph{local} high resolution 
is required.
\bibliographystyle{elsarticle-num}
\section*{\refname}

\appendix
\section{Detailed derivation of Pad\'e approximations}\label{sec:appendix}
\emph{The case $q=2$.}
The coefficients of this expansion are already known, however it is instructive
to proceed with their derivation in order to understand how it works.
We have to compute $3$ coefficients, $a_1$, $b_1$ and $a_2$, therefore we can
use only $3$ equations.
These equations are obtained by nullifying, respectively, the coefficients of the
terms $r^2$, $r^4$ and $r^6$ in the numerator $N_2(r)$ (lower-order powers of
$r^{2k}$).
By nullifying the coefficient of $r^2$, we get 
$$
-4 a_1^2 b_1 + 4 a_1 a_2 + a_1^2=0,
$$
from which $a_2=a_1\left(b_1-\frac{1}{4}\right)$.
By nullifying the coefficient of $r^4$, and replacing $a_2$ with the expression 
above, we get 
$$
a_1^2 (192 a_1 b_1-48 b_1-32 a_1+5)=0,
$$
which gives $b_1=\frac{5-32a_1}{48-192a_1}$.
If we now nullify the coefficient of $r^6$ and replace $a_2$ with 
$a_1\left(b_1-\frac{1}{4}\right)$ and $b_1$ with
$\frac{5-32a_1}{48-192a_1}$, we get the following equation 
$$
a_1^2 (8 a_1+1) (32 a_1-11) = 0.
$$
Clearly, $a_1=0$ is not acceptable, nor is $a_1= -\frac 1 8$. 
The only acceptable value is $a_1=\frac{11}{32}$.
As we mentioned before, equation~\eqref{eq:GPErho} cannot be satisfied exactly,
however, an \emph{a posteriori} evaluation reveals that the remaining 
coefficients
of $r^{2k}$ are smaller than $1.5 \times 10^{-4}$ and monotonically decreasing 
with $k$.

\begin{table}
{\footnotesize
$$
\begin{array}{c|cc}
& q=2 & q=3\\
\hline \\
a_1 & \Fra{11}{32} & 0.34003812123694735361\\
& & \\
b_1 & \Fra{5-32a_1}{48-192a_1} &
      \Fra{2304a_1^3+656a_1^2-421a_1-28}{7680a_1^2-1680a_1-330}  \\
& & \\
a_2 &  a_1\left(b_1-\Fra{1}{4}\right) &  a_1\left(b_1-\Fra{1}{4}\right) \\
& & \\
b_2 & & \Fra{768a_1b_1-120b_1-384a_1^2+8a_1+7}{4608a_1-1152} \\
& &\\
a_3 & & \Fra{a_1(192b_2-48b_1+16a_1+5)}{192} \\
\end{array}
$$}
\caption{Coefficients of Pad\'e approximations $\rho_2$ and $\rho_3$.}
\label{t:coeff23}
\end{table}

\emph{The case $q=3$.}
Since we have to compute $5$ coefficients we need $5$ equations, which are
obtained by imposing that the coefficients of the
terms $r^2$, $r^4$, $r^6$, $r^8$ and $r^{10}$ must be zero.
By nullifying the coefficient of $r^2$ we still get the same equation as for
$q=2$, $-4 a_1^2 b_1 + 4 a_1 a_2 + a_1^2=0$, from which
$a_2=a_1\left(b_1-\frac{1}{4}\right)$.
By nullifying the coefficient of $r^4$, and replacing $a_2$ with the expression 
above, we get 
$$
192 a_1 b_2 - 48 a_ 1 b_1 - 192 a_3 + 16 a_1^2 + 5 a_1 = 0,
$$
which is easy to solve for $a_3$ leading to
$$
a_3 = \frac{a_1(192b_2-48b_1+16a_1+5)}{192}.
$$
Now we collect terms in $r^6$ and impose its coefficient to be zero.
In this equation we replace $a_2$ and $a_3$ with the expressions derived 
above and get the equation
$$
4608 a_1  b_2 - 1152 b_2 - 768 a_1 b_1 + 120 b_1 + 384 a_1^2 - 8 a_1-7=0,
$$
which we solve for $b_2$:
$$
b_2=\Fra{768a_1b_1-120b_1-384a_1^2+8a_1+7}{4608a_1-1152}.
$$
Then we nullify the coefficient of $r^8$, substitute all previously found
$a_2$, $a_3$ and $b_2$, getting the equation 
$$
7680 a_1^2 b_1-1680 a_1 b_1-330 b_1-2304 a_1^3-656 a_1^2+421 a_1+28=0,
$$
which gives
$$
b_1=\Fra{2304a_1^3+656a_1^2-421a_1-28}{7680a_1^2-1680a_1-330}.
$$
Finally, we nullify the coefficient of $r^{10}$, substitute 
$a_2$, $a_3$, $b_2$ and $b_1$, and get the equation for $a_1$
$$
a_1^2\left(21233664a_1^5 - 9732096a_1^4 - 62464a_1^3 + 137856a_1^2 + 
62772a_1 - 1247 \right)=0.
$$
This equation must be solved numerically and leads to many real solutions.
However the only value that reproduces a physical behavior of $\rho_3(r)$ for $r
\to 0$ is $a_1=0.34003812123694735361$.
It is possible to compute the first derivative and verify that $\rho_3'(r)>0$ 
for all $r>0$.
In other words, $\rho_3$ is a physical, monotonically increasing, approximation
of the density due to a two-dimensional quantum vortex.
Again, the coefficients of $r^{2k}$ that are not zero are, indeed, smaller
than $4.0 \times 10^{-4}$ and monotonically decreasing with $k$.

\begin{table}
{\footnotesize
$$
\begin{array}{c|c}
& q=4\\
\hline \\
a_1 & 0.34010790700196714760\\
& \\
b_1 & \Fra{2304a_1^3+656a_1^2-421a_1-28}{7680a_1^2-1680a_1-330} \\
& \\
a_2 & a_1\left(b_1-\frac{1}{4}\right)\\
& \\
b_2 & \Fra{(737280 a_1^3+209920 a_1^2-134720 a_1-8960) b_1-364544 a_1^3+
70144 a_1^2+18256 a_1+393}{2457600 a_1^2-537600 a_1-105600} \\
& \\
a_3 & \Fra{a_1(192b_2-48b_1+16a_1+5)}{192}\\
& \\
b_3 & \Fra{(61440 a_1-9600) b_2+(-30720 a_1^2+640 a_1+560) b_1+8448 a_1^2-1056 a_1-21}
{368640 a_1-92160} \\
& \\
a_4 & \Fra{4608 a_1 b_3-1152 a_1 b_2+(384 a_1^2+120 a_1) b_1-128 a_1^2-7 a_1}{4608} \\
\end{array}
$$}
\caption{Coefficients of Pad\'e approximation $\rho_4$.}
\label{t:coeff4}
\end{table}

\emph{The case $q=4$.}
Now we have 7 coefficients to compute, therefore we need $7$ equations, i.e. we
need to nullify the coefficients of $r^{2k}$ for $k=1,\dots,7$.
By canceling the term $r^2$ and solving for $a_2$ we get the usual expression
$a_2=a_1\left(b_1-\frac{1}{4}\right)$.
By nullifying the term $r^4$, substituting $a_2$ and solving for $a_3$ we get
$a_3=\frac{a_1(192b_2-48b_1+16a_1+5)}{192}$, which is the same expression
obtained for $\rho_3$.
By canceling the term $r^6$, substituting $a_2$ and $a_3$ as found, and solving
for $a_4$ we get
$$
a_4=
\frac{4608 a_1 b_3-1152 a_1 b_2+(384 a_1^2+120 a_1) b_1-128 a_1^2-7 a_1}{4608}.
$$
By canceling the term $r^8$, substituting $a_2$, $a_3$, $a_4$ and solving
for $b_3$ we get
\small
$$
b_3=
\frac{(61440 a_1-9600) b_2+(-30720 a_1^2+640 a_1+560) b_1+8448 a_1^2-1056 a_1-21}
{368640 a_1-92160}.
$$
\normalsize
By canceling the term $r^{10}$, substituting $a_2$, $a_3$, $a_4$ and $b_3$, 
and solving for $b_2$ we get
\footnotesize
$$
b_2=
\frac{(737280 a_1^3+209920 a_1^2-134720 a_1-8960) b_1-364544 a_1^3+
70144 a_1^2+18256 a_1+393}{2457600 a_1^2-537600 a_1-105600}.
$$
\normalsize
By canceling the term $r^{12}$, substituting all known $a_j$, $b_3$, $b_2$, and 
solving for $b_1$ we get
\footnotesize
$$
b_1=
\frac{722731008 a_1^5-326467584 a_1^4-13427712 a_1^3+11551104 a_1^2+
834006 a_1-12183}{2972712960 a_1^5-1362493440 a_1^4-8744960 a_1^3+
19299840 a_1^2+8788080 a_1-174580}
$$
\normalsize
Finally by canceling the term $r^{14}$, substituting all $a_j$, $b_k$ previously
found, we get an equation for $a_1$
$$
\begin{aligned}
&1292033536819200 a_1^8-2530164294549504 a_1^7+1853440540016640 a_1^6-\\
&642522859438080 a_1^5+107808283328512 a_1^4-8028170208256 a_1^3+\\
&248539665024 a_1^2+1297120628 a_1+9325957=0.
\end{aligned}
$$
This equation has many real solutions, which can be determined numerically.
However, the only value that leads to a
physically acceptable $\rho_4(r)$ for $r \to 0$ is
$a_1=0.34010790700196714760$.
After computing all other coefficients and the first derivative, it is
straighforward to verify that $\rho_4'(r)>0$ for all $r>0$, i.e.
$\rho_4$ is a physical, monotonically increasing, approximation
of the density for a two-dimensional quantum vortex.
As observed before, the coefficients of $r^{2k}$ that are not zero are smaller
than $1.9 \times 10^{-11}$ and monotonically decreasing with $k$. 
\end{document}